\documentclass[10pt,twocolumn,twoside]{IEEEtran}              
\IEEEoverridecommandlockouts         
\usepackage{amsmath,amsfonts,amssymb,color}
\usepackage{epsfig}
\usepackage{psfrag}
\usepackage{algorithm}
\usepackage{algpseudocode}
\usepackage{array}
\usepackage{epstopdf}
\usepackage{cite}
\usepackage{footmisc}
\usepackage{mathtools}
\usepackage{graphicx}
\DeclarePairedDelimiter{\ceil}{\lceil}{\rceil}

\newtheorem{problem}{Problem}
\newtheorem{theorem}{Theorem}
\newtheorem{corollary}{Corollary}
\newtheorem{lemma}{Lemma}
\newtheorem{remark}{Remark}
\newtheorem{definition}{Definition}

\newtheorem{example}{Example}

\renewcommand{\theenumi}{(\alph{enumi}}
\algnewcommand{\LineComment}[1]{\State \(\triangleright\) #1}

\title{\LARGE \bf
On the Complexity and Approximability of Optimal Sensor Selection and Attack for Kalman Filtering
}
\begin{document}

\author{Lintao Ye, Nathaniel Woodford, Sandip Roy and Shreyas Sundaram% <-this % stops a space
\thanks{This research was supported by NSF grants CMMI-1635014 and CMMI-1635184. The work of the second author was also supported by the Purdue Military Research Initiative.  Lintao Ye, Nathaniel Woodford and Shreyas Sundaram are with the School of Electrical
and Computer Engineering at Purdue University, West Lafayette, IN 47906 USA. Email: \{ye159,nwoodfor,sundara2\}@purdue.edu.  Sandip Roy is with the School of Electrical Engineering and Computer Science at Washington State University, Pullman, WA 99164 USA.  Email: sroy@eecs.wsu.edu.
}
}
\maketitle
\thispagestyle{empty}
\pagestyle{empty}

%%%%%%%%%%%%%%%%%%%%%%%%%%%%%%%%%%%%%%%%%%%%%%%%%%%%%%%%%%%%%%%%%%%%%%%%%%%%%%%%
\begin{abstract}
Given a linear dynamical system affected by stochastic noise, we consider the problem of selecting an optimal set of sensors (at design-time) to minimize the trace of the steady state \textit{\textbf{a priori}} or \textit{\textbf{a posteriori}} error covariance of the Kalman filter, subject to certain selection budget constraints. We show the fundamental result that there is no polynomial-time constant-factor approximation algorithm for this problem. This contrasts with other classes of sensor selection problems studied in the literature, which typically pursue constant-factor approximations by leveraging greedy algorithms and submodularity (or supermodularity) of the cost function. Here, we provide a specific example showing that greedy algorithms can perform arbitrarily poorly for the problem of design-time sensor selection for Kalman filtering. We then study the problem of attacking (i.e., removing) a set of installed sensors, under predefined attack budget constraints, to maximize the trace of the steady state \textit{\textbf{a priori}} or \textit{\textbf{a posteriori}} error covariance of the Kalman filter. Again, we show that there is no  polynomial-time constant-factor approximation algorithm for this problem, and show specifically that greedy algorithms can perform arbitrarily poorly.
\end{abstract}

%%%%%%%%%%%%%%%%%%%%%%%%%%%%%%%%%%%%%%%%%%%%%%%%%%%%%%%%%%%%%%%%%%%%%%%%%%%%%%%%
\section{Introduction}
In large-scale control system design, the number of sensors or actuators that can be installed is typically limited by a design budget constraint. Moreover, system designers often need to select among a set of possible sensors and actuators, with varying qualities and costs. Consequently, a key problem is to select an appropriate set of sensors or actuators in order to achieve certain objectives. This problem has recently received much attention from researchers (e.g., \cite{van2001review,olshevsky2014minimal,liu2003problem,joshi2009sensor,tzoumas2016minimal,summers2016submodularity,olshevsky2018non,lin2013design,yoo2002np}). One specific instance of this problem arises in the context of linear Gauss-Markov systems, where the corresponding Kalman filter (with the chosen sensors) is used to estimate the states of the systems (e.g., \cite{mo2011sensor}, \cite{chamon2020approximate}). The problem then becomes how to select sensors dynamically (at run-time) or statically (at design-time) to minimize certain metrics of the corresponding Kalman filter. The former scenario is known as the sensor scheduling problem, where different sets of sensors can be chosen at different time steps (e.g., \cite{vitus2012efficient,huber2012optimal,jawaid2015submodularity}). The latter scenario is known as the design-time sensor selection problem, where the set of sensors is chosen a priori and is not allowed to change over time (e.g., \cite{dhingra2014admm,yang2015deterministic,ye2018optimal}). 

Since these problems are NP-hard in general (e.g., \cite{zhang2017sensor}), approximation algorithms that provide solutions within a certain factor of the optimal are then proposed to tackle them. Among these approximation algorithms, greedy algorithms have been widely used (e.g, \cite{krause2008near}, \cite{shamaiah2010greedy}), since such algorithms have provable performance guarantees if the cost function is submodular or supermodular (e.g., \cite{nemhauser1978analysis}, \cite{das2011submodular}).

Additionally, in many applications, the sensors that have been selected and installed on the system are susceptible to a variety of potential attacks. For instance, an adversary (attacker) can inject false data to corrupt the state estimation, which is known as the false data injection attack (e.g., \cite{liu2011false,mo2016performance,yang2014false}). Another type of attack is the Denial-of-Service (DoS) attack, where an attacker tries to diminish or eliminate a network's capacity to achieve its expected objective \cite{wood2002denial}, including, for example, wireless jamming (e.g., \cite{xu2006jamming}, \cite{zhang2016optimal}) and memory exhaustion through flooding (e.g., \cite{schuba1997analysis}). One class of DoS attacks corresponds to removing a set of installed sensors from the system, i.e., the measurements of the attacked sensors are not used. This was also studied in \cite{laszka2015resilient} and \cite{tzoumas2017resilient}, and will be the type of attack we consider in this work.

In this paper, we consider the sensor selection problem and the sensor attack problem for Kalman filtering of discrete-time linear dynamical systems. First, we study the problem of choosing a subset of sensors to install (under given selection budget constraints) to minimize the trace of either the steady state {\it a priori} or {\it a posteriori} error covariance of the Kalman filter.  We refer to these problems as the priori and posteriori {\it Kalman Filtering Sensor Selection (KFSS)} problems, respectively. Second, we study the problem of attacking the installed sensors (by removing a subset of them, under given attack budget constraints) to maximize the trace of either the steady state {\it a priori} or {\it a posteriori} error covariance of the Kalman filter associated with the surviving sensors. These problems are denoted as the priori and posteriori {\it Kalman Filtering Sensor Attack (KFSA)} problems, respectively. 

\subsection*{Contributions}
Our contributions are as follows. First, we show that for the priori and posteriori KFSS problems, {\it there are no  polynomial-time constant-factor approximation algorithms for these problems} (unless P $=$ NP) {\it even} for the special case when the system is stable and all sensors have the same cost. In other words, there are no polynomial-time algorithms that can find a sensor selection that is always guaranteed to yield a mean square estimation error (MSEE) that is within any constant finite factor of the MSEE for the optimal selection.  More importantly, our result stands in stark contrast to other sensor selection problems studied in the literature, which leveraged submodularity of their associated cost functions to provide greedy algorithms with constant-factor approximation ratios (e.g., \cite{krause2008near}, \cite{tzoumas2016sensor}). Second, we show that the same results hold for the priori and posteriori KFSA problems. Our inapproximability results directly imply that greedy algorithms cannot provide constant-factor guarantees for our problems.  Our third contribution is to explicitly show how greedy algorithms can provide arbitrarily poor performance even for very small instances (with only three states) of the priori and posteriori KFSS (resp., KFSA) problems. 

A portion of the results pertaining to only the priori KFSS problem appears in the conference paper \cite{ye2018complexity}.

\subsection*{Related work}
The authors in \cite{chamon2020approximate} and \cite{tzoumas2016sensor} studied the design-time sensor selection problem for discrete-time linear time-varying systems over a finite time horizon. The objective is to minimize the estimation error with a cardinality constraint on the chosen sensors (or alternatively, minimize the number of chosen sensors while guaranteeing a certain level of performance in terms of the estimation error). The authors analyzed the performance of greedy algorithms for this problem. However, their results cannot be directly applied to the problems that we consider here, since we aim to optimize the steady state estimation error.  

The papers \cite{yang2015deterministic} and \cite{zhang2017sensor} considered the same design-time sensor selection as the one we consider here. In \cite{yang2015deterministic}, the authors expressed the problem as a semidefinite program (SDP). However, they did not provide theoretical guarantees on the performance of the proposed algorithm. The paper \cite{zhang2017sensor} showed that the problem is NP-hard and provided upper bounds on the performance of any algorithm for the problem; these upper bounds were functions of the system matrices. Although \cite{zhang2017sensor} showed via simulations that greedy algorithms performed well for several randomly generated systems, the question of whether such algorithms (or other polynomial-time algorithms) could provide constant-factor approximation ratios for the problem was left open. We resolve this question in this paper by showing that there does not exist any  polynomial-time constant-factor approximation algorithm for this problem.

In \cite{laszka2015resilient}, the authors studied the problem of attacking a given observation selection in Gaussian process regression \cite{das2008algorithms} to maximize the {\it posteriori} variance of the predictor variable. It was shown that this problem is NP-hard. Moreover, they also gave an instance of this problem such that a greedy algorithm for finding an optimal attack will perform arbitrarily poorly. In \cite{lu2016sparse}, the authors considered the scenario where the attacker can target a different set of sensors at each time step to maximize certain metrics of the error covariance of the Kalman filter at the final time step. Some suboptimal algorithms were provided with simulation results. Different from \cite{laszka2015resilient} and \cite{lu2016sparse}, we study the problem where the attacker removes a set of installed sensors to maximize the trace of the steady state error covariance of the Kalman filter associated with the surviving sensors, and provide fundamental limitations on achievable performance by any possible algorithm for this problem.

\subsection*{Notation and Terminology}
The sets of integers and real numbers are denoted as $\mathbb{Z}$ and $\mathbb{R}$, respectively. The set of integers that are greater than (resp., greater than or equal to) $a\in\mathbb{R}$ is denoted as $\mathbb{Z}_{>a}$ (resp., $\mathbb{Z}_{\ge a}$). Similarly, we use the notations $\mathbb{R}_{>a}$ and $\mathbb{R}_{\ge a}$. For any $x\in\mathbb{R}$, let $\ceil{x}$ be the least integer greater than or equal to $x$. For a square matrix $P\in\mathbb{R}^{n\times n}$, let $P^T$, $\textrm{rank}(P)$, $\textrm{rowspace}(P)$ and $\textrm{trace}(P)$ be its transpose, rank, rowspace and trace, respectively.  We use $P_{ij}$ (or $(P)_{ij}$) to denote the element in the  $i$th row and $j$th column of $P$. A diagonal matrix $P\in\mathbb{R}^{n\times n}$ is denoted as $\text{diag}(P_{11},\dots,P_{nn})$. The set of $n$ by $n$ positive definite (resp., positive semi-definite) matrices is denoted as $\mathbb{S}_{++}^n$ (resp., $\mathbb{S}_{+}^n$). The identity matrix with dimension $n\times n$ is denoted as $I_{n}$. The zero matrix with dimension $m\times n$ is denoted as $\mathbf{0}_{m\times n}$. In a matrix, let $*$ denote elements of the matrix that are of no interest. For a vector $v$, let $v_i$ be the $i$th element of $v$ and define the support of $v$ to be $\textrm{supp}(v)=\{i:v_i\neq 0\}$.  The Euclidean norm of $v$ is denoted by $\lVert v\lVert_2$. Define $\mathbf{e}_i$ to be a row vector where the $i$th element is $1$ and all the other elements are zero; the dimension of the vector can be inferred from the context. Define $\mathbf{1}_n$ to be a column vector of dimension $n$ with all the elements equal to 1. The set of $0$-$1$ indicator vectors of dimension $n$ is denoted as $\{0,1\}^n$. For a random variable $\omega$, let $\mathbb{E}[\omega]$ be its expectation. For a set $\mathcal{A}$, let $|\mathcal{A}|$ be its cardinality.

\section{Problem Formulation} \label{sec:problem formulation}
Consider the discrete-time linear system
\begin{equation}
x[k+1] = Ax[k] + w[k],
\label{eqn:system dynamics}
\end{equation}
where $x[k]\in\mathbb{R}^n$ is the system state, $w[k]\in\mathbb{R}^n$ is a zero-mean white Gaussian noise process with $\mathbb{E}[w[k](w[k])^T]=W$, and $A\in\mathbb{R}^{n\times n}$ is the system dynamics matrix. The initial state $x[0]$ is assumed to be a Gaussian random vector with mean $\bar{x}_0\in\mathbb{R}^n$ and covariance $\Pi_0\in\mathbb{S}_{+}^n$.

Consider a set $\mathcal{Q}$ that contains $q$ sensors. Each sensor $i\in\mathcal{Q}$ provides a measurement of the system of the form
\begin{equation}
\label{eqn:single sensor measurement}
y_i[k]=C_ix[k]+v_i[k],
\end{equation}
where $C_i\in\mathbb{R}^{s_i\times n}$ is the measurement matrix for sensor $i$, and $v_i[k]\in\mathbb{R}^{s_i}$ is a zero-mean white Gaussian noise process. We further define $y[k]\triangleq\big[(y_1[k])^T\ \cdots \ (y_q[k])^T\big]^T$, $C\triangleq\big[C_1^T\ \cdots\ C_q^T\big]^T$ and $v[k]\triangleq\big[(v_1[k])^T\ \cdots \ (v_q[k])^T\big]^T$. Thus, the output provided by all sensors together is given by
\begin{equation}
\label{eqn:all sensors measurements}
y[k]=C x[k]+v[k],
\end{equation}
where $C\in\mathbb{R}^{s\times n}$ and $s=\sum_{i=1}^q s_i$. Denote $\mathbb{E}[v[k](v[k])^T]=V$. We assume that the system noise and the measurement noise are uncorrelated, i.e., $\mathbb{E}[v[k](w[j])^T]=\mathbf{0}_{s\times n}$, $\forall k, j\in\mathbb{Z}_{\ge0}$, and $x[0]$ is independent of $w[k]$ and $v[k]$, $\forall k\in\mathbb{Z}_{\ge0}$.

\subsection{The Sensor Selection Problem}

Consider the scenario where there are no sensors initially deployed on the system. Instead, the system designer must select a subset of sensors from $\mathcal{Q}$ to install. Each sensor $i\in\mathcal{Q}$ has a cost $b_{i}\in\mathbb{R}_{\ge0}$; define the cost vector $b\triangleq \left[\begin{matrix}b_1 & \cdots & b_q\end{matrix}\right]^T$. The designer has a budget $B\in\mathbb{R}_{\ge0}$ that can be spent on choosing sensors from $\mathcal{Q}$.

After a set of sensors is selected and installed, the Kalman filter is applied to provide an optimal estimate of the states using the measurements from the installed sensors (in the sense of minimizing the MSEE). Define $\mu\in\{0,1\}^q$ to be the indicator vector of the selected sensors, where $\mu_i=1$ if and only if sensor $i\in\mathcal{Q}$ is installed. Let $C(\mu)$ denote the measurement matrix of the installed sensors indicated by $\mu$, i.e., $C(\mu)\triangleq\left[\begin{matrix} C_{i_1}^T &  \cdots & C_{i_p}^T \end{matrix}\right]^T$, where $\textrm{supp}(\mu)=\{i_1,\dots,i_p\}$. Similarly, let $V(\mu)$ be the measurement noise covariance matrix of the installed sensors, i.e., $V(\mu)=\mathbb{E}[\tilde{v}[k](\tilde{v}[k])^T]$, where $\tilde{v}[k]=\big[(v_{i_1}[k])^T\ \cdots \ (v_{i_p}[k])^T\big]^T$. Let $\Sigma_{k|k-1}(\mu)$ (resp., $\Sigma_{k|k}(\mu)$) denote the {\it a priori} (resp., {\it a posteriori}) error covariance matrix of the Kalman filter at time step $k$, when the sensors indicated by $\mu$ are installed. Take the initial covariance $\Sigma_{0|-1}(\mu)=\Pi_0$, $\forall\mu$.  We will use the following result \cite{anderson1979optimal}.

\begin{lemma}
Suppose that the pair $(A,W^{\frac{1}{2}})$ is stabilizable. For a given indicator vector $\mu$, $\Sigma_{k|k-1}(\mu)$ (resp., $\Sigma_{k|k}(\mu)$) will converge to a finite limit $\Sigma(\mu)$ (resp., $\Sigma^*(
\mu)$), which does not depend on the initial covariance $\Pi_0$, as $k\to\infty$ if and only if the pair $(A,C(\mu))$ is detectable.
\label{lemma:Anderson optimal filtering}
\end{lemma}

The limit $\Sigma(\mu)$ satisfies the \textit{discrete algebraic Riccati equation (DARE)} \cite{anderson1979optimal}:
\begin{multline}
\label{eqn:DARE}
\Sigma(\mu)=A\Sigma(\mu)A^T+W - \\
A\Sigma(\mu)C(\mu)^T\big(C(\mu)\Sigma(\mu)C(\mu)^T+V(\mu)\big)^{-1}C(\mu)\Sigma(\mu)A^T.
\end{multline}
The limits $\Sigma(\mu)$ and $\Sigma^*(\mu)$ are coupled as
\begin{equation}
\label{eqn:couple priori and posteriori 1}
\Sigma(\mu)=A\Sigma^*(\mu)A^T+W.
\end{equation}
The limit $\Sigma^*(\mu)$ of the a {\it posteriori} error covariance matrix satisfies the following equation \cite{catlin1989estimation}:
\begin{multline}
\label{eqn:postDARE}
\Sigma^*(\mu)=\Sigma(\mu)-\\
\Sigma(\mu)C(\mu)^T(C(\mu)\Sigma(\mu)C(\mu)^T+V(\mu))^{-1}C(\mu)\Sigma(\mu).
\end{multline}
Note that we can either obtain $\Sigma^*(\mu)$ from $\Sigma(\mu)$ using Eq. \eqref{eqn:postDARE} or by substituting Eq. \eqref{eqn:couple priori and posteriori 1} into Eq. \eqref{eqn:postDARE} and solving for $\Sigma^*(\mu)$. The inverses in Eq. \eqref{eqn:DARE} and Eq. \eqref{eqn:postDARE} are interpreted as pseudo-inverses if the arguments are not invertible. 

For the case when the pair $(A,C(\mu))$ is not detectable, we define $\Sigma(\mu)=+\infty$ and $\Sigma^*(\mu)=+\infty$. Moreover, for any sensor selection $\mu$, we note from Lemma $\ref{lemma:Anderson optimal filtering}$ that the limit $\Sigma(\mu)$ (resp., $\Sigma^*(\mu)$), if it exists, does not depend on $\bar{x}_0$ or $\Pi_0$. Thus, we can assume without loss of generality that $\bar{x}_0=\mathbf{0}$ and $\Pi_0=I_n$ in the sequel.  The priori and posteriori Kalman Filtering Sensor Selection (KFSS) problems are then defined as follows.
\begin{problem}
\label{pro:priori KFSS problem}
\label{pro:posteriori KFSS problem}
(Priori and Posteriori KFSS Problems). Given a system dynamics matrix $A\in\mathbb{R}^{n\times n}$, a measurement matrix $C\in\mathbb{R}^{s\times n}$ containing all of the individual sensor measurement matrices, a system noise covariance matrix $W\in\mathbb{S}_+^n$, a sensor noise covariance matrix $V\in\mathbb{S}_+^{s}$, a cost vector $b\in\mathbb{R}^q_{\ge0}$ and a budget $B\in\mathbb{R}_{\ge0}$, the priori Kalman filtering sensor selection problem is to find the sensor selection $\mu$, i.e., the indicator vector $\mu$ of the selected sensors, that solves
\begin{equation*}
\begin{split}
&\mathop{\min}_{\mu\in\{0,1\}^q}\ \text{trace}(\Sigma(\mu))\\
&s.t.\ b^T \mu\le B
\end{split}
\end{equation*}
where $\Sigma(\mu)$ is given by Eq. $(\ref{eqn:DARE})$ if the pair $(A,C(\mu))$ is detectable, and $\Sigma(\mu)=+\infty$ if otherwise. Similarly, the posteriori Kalman filtering sensor selection problem is to find the sensor selection $\mu$ that solves
\begin{equation*}
\begin{split}
&\mathop{\min}_{\mu\in\{0,1\}^q}\ \text{trace}(\Sigma^*(\mu))\\
&s.t.\ b^T \mu\le B
\end{split}
\end{equation*}
where $\Sigma^*(\mu)$ is given by Eq. $(\ref{eqn:postDARE})$ if the pair $(A,C(\mu))$ is detectable, and $\Sigma^*(\mu)=+\infty$ if otherwise.
\end{problem}

\subsection{The Sensor Attack Problem}
Now consider the scenario where the set $\mathcal{Q}$ of sensors has already been installed on the system. An adversary desires to attack a subset of sensors (i.e., remove a subset of sensors from the system), where each sensor $i\in\mathcal{Q}$ has an attack cost $\omega_i\in\mathbb{R}_{\ge0}$; define the cost vector $\omega\triangleq \left[\begin{matrix}\omega_1 & \cdots & \omega_q\end{matrix}\right]^T$. We assume that the adversary has a budget $\Omega\in\mathbb{R}_{\ge0}$, which is the total cost that can be spent on removing sensors from $\mathcal{Q}$.

After a subset of sensors are attacked (i.e., removed), the Kalman filter is then applied to estimate the states using the measurements from the surviving sensors. We define a vector $\nu\in\{0,1\}^q$ to be the indicator vector of the attacked sensors, where $\nu_i=1$ if and only if sensor $i\in\mathcal{Q}$ is attacked. Hence, the set of sensors that survive is $\mathcal{Q}\setminus\textrm{supp}(\nu)$. Define $v^c\in\{0,1\}^q$ to be the vector such that $\textrm{supp}(\nu^c)=\mathcal{Q}\setminus\textrm{supp}(\nu)$, i.e., $\nu^c_i=1$ if and only if sensor $i\in\mathcal{Q}$ survives. Similarly to the sensor selection problem, we let $C(\nu^c)$ and $V(\nu^c)$ denote the measurement matrix and the measurement noise covariance matrix, respectively, corresponding to $\nu^c$. Furthermore, let $\Sigma_{k|k-1}(\nu^c)$ and $\Sigma_{k|k}(\nu^c)$ denote the a {\it priori} error covariance matrix and the a {\it posteriori} error covariance matrix of the Kalman filter at time step $k$, respectively. Denote $\displaystyle\mathop{\lim}_{k\to\infty}\Sigma_{k|k-1}(\nu^c)=\Sigma(\nu^c)$ and $\displaystyle\mathop{\lim}_{k\to\infty}\Sigma_{k|k}(\nu^c)=\Sigma^*(\nu^c)$ if the limits exist, according to Lemma $\ref{lemma:Anderson optimal filtering}$.  Note that Eq. \eqref{eqn:DARE}-\eqref{eqn:postDARE} also hold if we substitute $\mu$ with $\nu^c$.

For the case when the pair $(A,C(\nu^c))$ is not detectable, we define $\Sigma(\nu^c)=+\infty$ and $\Sigma^*(\nu^c)=+\infty$. Recall that we have assumed without loss of generality that $\bar{x}_0=\mathbf{0}$ and $\Pi_0=I_n$.  The priori and posteriori Kalman Filtering Sensor Attack (KFSA) problems are defined as follows.
\begin{problem}
\label{pro:priori KFSA problem}
\label{pro:posteriori KFSA problem}
(Priori and Posteriori KFSA Problems). Given a system dynamics matrix $A\in\mathbb{R}^{n\times n}$, a measurement matrix $C\in\mathbb{R}^{s\times n}$, a system noise covariance matrix $W\in\mathbb{S}_+^n$, a sensor noise covariance matrix $V\in\mathbb{S}_+^{s}$, a cost vector $\omega\in\mathbb{R}^q_{\ge0}$ and a budget $\Omega\in\mathbb{R}_{\ge0}$, the priori Kalman filtering sensor attack problem is to find the sensor attack $\nu$, i.e., the indicator vector $\nu$ of the attacked sensors, that solves
\begin{equation*}
\begin{split}
&\mathop{\max}_{\nu\in\{0,1\}^q}\ \text{trace}(\Sigma(\nu^c))\\
&s.t.\ \omega^T \nu\le \Omega
\end{split}
\end{equation*}
where $\Sigma(\nu^c)$ is given by Eq. $(\ref{eqn:DARE})$ if the pair $(A,C(\nu^c))$ is detectable, and $\Sigma(\nu^c)=+\infty$ if otherwise. Similarly, the posteriori Kalman filtering sensor attack problem is to find the sensor attack $\nu$ that solves
\begin{equation*}
\begin{split}
&\mathop{\max}_{\nu\in\{0,1\}^q}\ \text{trace}(\Sigma^*(\nu^c))\\
&s.t.\ \omega^T \nu\le \Omega
\end{split}
\end{equation*}
where $\Sigma^*(\nu^c)$ is given by Eq. $(\ref{eqn:postDARE})$ if the pair $(A,C(\nu^c))$ is detectable, and $\Sigma^*(\nu^c)=+\infty$ if otherwise.
\end{problem}

Note that although we focus on the optimal sensor selection and attack problems for Kalman filtering, due to the duality between the Kalman filter and the Linear Quadratic Regulator (LQR) \cite{anderson2007optimal}, all of the analysis in this paper will also apply if the priori KFSS and KFSA problems are rephrased as optimal actuator selection and attack problems for LQR, respectively. We omit the details of the rephrasing in the interest of space. 

\begin{remark}
\label{remark:special instances}
Our goal in this paper is to show that for the priori and posteriori KFSS problems and the priori and posteriori KFSA problems, the optimal solutions cannot be approximated within any constant factor in polynomial time. To do this, it is sufficient for us to consider the special case when $C_i\in\mathbb{R}^{1\times n}$, $\forall i\in\{1,\dots,q\}$, i.e., each sensor provides a scalar measurement. Moreover, the sensor selection cost vector and the sensor attack cost vector are considered to be $b=\begin{bmatrix}1 & \cdots & 1\end{bmatrix}^T$ and $\omega=\begin{bmatrix}1 & \cdots & 1\end{bmatrix}^T$, respectively, i.e., the selection cost and the attack cost of each sensor are both equal to $1$. By showing that the problems are inapproximable even for these special subclasses, we obtain that the general versions of the problems are inapproximable as well.
\end{remark}

\section{Inapproximability of the KFSS and KFSA problems}\label{sec:complexity analysis KFSS and KFSA}
In this section, we analyze the approximability of the KFSS and KFSA problems. We will start with a brief overview of some relevant concepts from the field of computational complexity, and then provide some preliminary lemmas that we will use in proving our results. That will lead into our characterizations of the complexity of KFSS and KFSA.

\subsection{Review of Complexity Theory}
We first review the following fundamental concepts from complexity theory \cite{garey1979computers}. 
\begin{definition}
A {\it polynomial-time algorithm} for a problem is an algorithm that returns a solution to the problem in a polynomial (in the size of the problem) number of computations. 
\end{definition}
\begin{definition}
A {\it decision problem} is a problem whose answer is ``yes'' or ``no''. The set P contains those decision problems that can be solved by a polynomial-time algorithm. The set NP contains those decision problems whose ``yes'' answers can be verified using a polynomial-time algorithm.
\end{definition}
\begin{definition}
An {\it optimization problem} is a problem whose objective is to maximize or minimize a certain quantity, possibly subject to constraints.
\end{definition}
\begin{definition}
A problem $\mathcal{P}_1$ is {\it NP-complete} if (a) $\mathcal{P}_1\in$ NP and (b) for any problem $\mathcal{P}_2$ in NP, there exists a polynomial-time algorithm that converts (or ``reduces'') any instance of $\mathcal{P}_2$ to an instance of $\mathcal{P}_1$ such that the answer to the constructed instance of $\mathcal{P}_1$ provides the answer to the instance of $\mathcal{P}_2$. $\mathcal{P}_1$ is {\it NP-hard} if it satisfies (b), but not necessarily (a).
\end{definition}

The above definition indicates that if one had a polynomial-time algorithm for an NP-complete (or NP-hard) problem, then one could solve {\it every} problem in NP in polynomial time.  Specifically, suppose we had a polynomial-time algorithm to solve an NP-hard problem $\mathcal{P}_1$.  Then, given any problem $\mathcal{P}_2$ in NP, one could first reduce any instance of $\mathcal{P}_2$ to an instance of $\mathcal{P}_1$ in polynomial time (such that the answer to the constructed instance of $\mathcal{P}_1$ provides the answer to the given instance of $\mathcal{P}_2$), and then use the polynomial-time algorithm for $\mathcal{P}_1$ to obtain the answer to $\mathcal{P}_2$.  

The above discussion also reveals that to show that a given problem $\mathcal{P}_1$ is NP-hard, one simply needs to show that any instance of some other {\it NP-hard (or NP-complete)} problem $\mathcal{P}_2$ can be reduced to an instance of $\mathcal{P}_1$ in polynomial time (in such a way that the answer to the constructed instance of $\mathcal{P}_1$ provides the answer to the given instance of $\mathcal{P}_2$). For then, an algorithm for $\mathcal{P}_1$ can be used to solve $\mathcal{P}_2$, and hence, to solve all problems in NP (by NP-hardness of $\mathcal{P}_2$).

The following is a fundamental result in computational complexity theory \cite{garey1979computers}.

\begin{lemma}
\label{lemma:no polynomial algo}
If P $\neq$ NP, there is no polynomial-time algorithm for any NP-complete (or NP-hard) problem.
\end{lemma}

For NP-hard optimization problems, polynomial-time approximation algorithms are of particular interest. A constant-factor approximation algorithm is defined as follows.
\begin{definition}
\label{def:constant approx alg}
A {\it constant-factor approximation algorithm} for an optimization problem is an algorithm that always returns a solution within a constant (system-independent) factor of the optimal solution.
\end{definition}

We will discuss the notion of a constant-factor approximation algorithm in greater depth later in this section. As described in the Introduction, the KFSS problem was shown to be NP-hard in \cite{zhang2017sensor} for two classes of systems and sensor costs.  First, when the $A$ matrix is unstable, it was shown in \cite{olshevsky2014minimal} that the problem of selecting a subset of sensors to make the system detectable is NP-hard, which implies that KFSS is NP-hard using Lemma $\ref{lemma:Anderson optimal filtering}$ as shown in  \cite{zhang2017sensor}. Second, when the $A$ matrix is stable (so that all sensor selections cause the system to be detectable), \cite{zhang2017sensor} showed that when the sensor selection costs can be arbitrary, the knapsack problem can be encoded as a special case of the KFSS problem, thereby again showing NP-hardness of the latter problem.

In this paper, we will show that the hardness of KFSS (resp., KFSA) does not solely come from selecting (resp., attacking) sensors to make the system detectable (resp., undetectable) or the sensor selection (resp., attack) costs. To do this, we will show a stronger result that there is no  polynomial-time constant-factor approximation algorithm for KFSS (resp., KFSA) {\it even} when the corresponding system dynamics matrix $A$ is {\it stable} (which guarantees the detectability of the system), and all the sensors have the same selection (resp., attack) cost.  Specifically, we consider a known NP-complete problem, and show how to reduce it to certain instances of KFSS (resp., KFSA) with stable $A$ matrices in polynomial time such that hypothetical  polynomial-time constant-factor approximation algorithms for the latter problems can be used to solve the known NP-complete problem. Since we know from Lemma $\ref{lemma:no polynomial algo}$ that if P $\neq$ NP, there does not exist a polynomial-time algorithm for any NP-complete problem, we conclude that if P $\neq$ NP, there is no  polynomial-time constant-factor approximation algorithm for KFSS (resp., KFSA), which directly implies that the KFSS (resp., KFSA) problem is NP-hard even under the extra conditions described above. We emphasize that our results do not imply that there is no  polynomial-time constant-factor approximation algorithm for {\it specific} instances of KFSS (resp., KFSA). Rather, the result is that we cannot have such an algorithm for {\it all} instances of KFSS (resp., KFSA).

\subsection{Preliminary Results}
The following results characterize properties of the KFSS and KFSA instances that we will consider when proving the inapproximability of the KFSS and KFSA problems.  The proofs are provided in Appendix A. 
\begin{lemma}
Consider a discrete-time linear system defined in $\eqref{eqn:system dynamics}$ and $\eqref{eqn:all sensors measurements}$. Suppose the system dynamics matrix is of the form $A=\textrm{diag}(\lambda_1,\dots,\lambda_n)$ with $0\le |\lambda_i|<1$, $\forall i \in \{1,\dots,n\}$, the system noise covariance matrix $W\in\mathbb{S}_{+}^n$ is diagonal, and the sensor noise covariance matrix $V\in\mathbb{S}^q_{+}$. Then, the following hold for all sensor selections $\mu$.

(a) For all $i\in\{1,\dots,n\}$, $(\Sigma(\mu))_{ii}$ and $(\Sigma^*(\mu))_{ii}$ satisfy 
\begin{equation}
W_{ii}\le (\Sigma(\mu))_{ii}\le \frac{W_{ii}}{1-\lambda_i^2},
\label{eqn:bounds for priori sigma}
\end{equation}
and
\begin{equation}
0\le(\Sigma^*(\mu))_{ii}\le\frac{W_{ii}}{1-\lambda_i^2},
\label{eqn:bounds for posteriori sigma}
\end{equation}
respectively.

(b) If $\exists i\in\{1,\dots,n\}$ such that $W_{ii}\neq0$ and the $i$th column of $C$ is zero, then $(\Sigma(\mu))_{ii}=(\Sigma^*(\mu))_{ii}=\frac{W_{ii}}{1-\lambda_i^2}$.

(c) If $V=\mathbf{0}_{q\times q}$ and $\exists i\in\{1,\dots,n\}$ such that $\mathbf{e}_i\in\textrm{rowspace}(C(\mu))$, then $(\Sigma(\mu))_{ii}=W_{ii}$ and $(\Sigma^*(\mu))_{ii}=0$.
\label{Lemma:minimum trace of sigma}
\end{lemma}

\begin{lemma}
Consider a discrete-time linear system defined in Eq. $(\ref{eqn:system dynamics})$ and Eq. $(\ref{eqn:all sensors measurements})$. Suppose the system dynamics matrix is of the form $A=\textrm{diag}(\lambda_1,0,\dots,0)\in\mathbb{R}^{n\times n}$, where $0<|\lambda_1|<1$, and the system noise covariance matrix is $W=I_n$. 

(a) Suppose the measurement matrix is of the form $C=\begin{bmatrix}1 & \mathbf{\gamma}\end{bmatrix}$ with sensor noise variance $V=\sigma_v^2$, where $\gamma\in\mathbb{R}^{1\times(n-1)}$ and $\sigma_v\in\mathbb{R}_{\ge0}$. Then, the MSEE of state $1$, denoted as $\Sigma_{11}$, satisfies 
\begin{equation}
\Sigma_{11}=\frac{1+\alpha^2\lambda_1^2-\alpha^2+\sqrt[]{(\alpha^2-\alpha^2\lambda_1^2-1)^2+4\alpha^2}}{2},
\label{eqn: mean square estimation error of state 1}
\end{equation}
where $\alpha^2=\lVert\gamma\rVert_2^2+\sigma_v^2$. 

(b) Suppose the measurement matrix is of the form $C=\left[\begin{matrix}\mathbf{1}_{n-1} & \rho I_{n-1}\end{matrix}\right]$ with sensor noise covariance $V=\mathbf{0}_{(n-1)\times(n-1)}$, where $\rho\in\mathbb{R}$. Then, the MSEE of state $1$, denoted as $\Sigma'_{11}$, satisfies
\begin{equation}
\Sigma'_{11}=\frac{\lambda_1^2\rho^2+n'-\rho^2+\sqrt{(\rho^2-\lambda_1^2\rho^2-n')^2+4n'\rho^2}}{2n'},
\label{eqn:mean square estimation error of state 1 with identity noise states}
\end{equation}
where $n'=n-1$. 

Moreover, if we view $\Sigma_{11}$ and $\Sigma'_{11}$ as functions of $\alpha^2$ and $\rho^2$, denoted as $\Sigma_{11}(\alpha^2)$ and $\Sigma'_{11}(\rho^2)$, respectively, then $\Sigma_{11}(\alpha^2)$ and $\Sigma'_{11}(\rho^2)$ are strictly increasing functions of $\alpha^2\in\mathbb{R}_{\ge0}$ and $\rho^2\in\mathbb{R}_{\ge0}$, with ${\lim}_{\alpha\to\infty}\Sigma_{11}(\alpha^2)=\frac{1}{1-\lambda_1^2}$ and ${\lim}_{\rho\to\infty}\Sigma'_{11}(\rho^2)=\frac{1}{1-\lambda_1^2}$, respectively.
\label{lemma:estimation error of state 1}
\end{lemma}

\subsection{Inapproximability of the KFSS Problem}
In this section, we characterize the achievable performance of algorithms for the priori and posteriori KFSS problems. For any given algorithm $\mathcal{A}$ (resp., $\mathcal{A}'$) of the priori (resp., posteriori) KFSS problem, we define the following ratios:
\begin{equation}
r_{\mathcal{A}}(\Sigma)\triangleq \frac{\textrm{trace}(\Sigma_{\mathcal{A}})}{\textrm{trace}(\Sigma_{opt})},
\label{eqn:performance ratio of algorithm}
\end{equation}
and
\begin{equation}
r_{\mathcal{A'}}(\Sigma^*)\triangleq \frac{\textrm{trace}(\Sigma^*_{\mathcal{A'}})}{\textrm{trace}(\Sigma^*_{opt})},
\label{eqn:performance ratio of algorithm post}
\end{equation}
where $\Sigma_{opt}$ (resp., $\Sigma^*_{opt}$) is the optimal solution to the priori (resp., posteriori) KFSS problem and $\Sigma_{\mathcal{A}}$ (resp., $\Sigma^*_{\mathcal{A}'}$) is the solution to the priori (resp., posteriori) KFSS problem given by algorithm $\mathcal{A}$ (resp., $\mathcal{A}'$). 

In \cite{zhang2017sensor}, the authors showed that there is an upper bound for $r_{\mathcal{A}}(\Sigma)$ (resp., $r_\mathcal{A'}(\Sigma^*)$) for any sensor selection algorithm $\mathcal{A}$ (resp., $\mathcal{A}'$), in terms of the system matrices.  However, the question of whether it is possible to find an algorithm $\mathcal{A}$ (resp., $\mathcal{A}'$) that is guaranteed to provide an approximation ratio $r_{\mathcal{A}}(\Sigma)$ (resp., $r_{\mathcal{A}'}(\Sigma^*)$) that is {\it independent} of the system parameters has remained open up to this point. In particular, it is desirable to find polynomial-time {\it constant-factor} approximation algorithms for the priori (resp., posteriori) KFSS problem, where the ratio $r_{\mathcal{A}}(\Sigma)$ (resp., $r_{\mathcal{A}'}(\Sigma^*)$) is upper-bounded by some (system-independent) constant.\footnote{Polynomial-time constant-factor approximation algorithms have been widely studied for NP-hard problems, e.g.,\cite{garey1979computers}.} We provide a negative result by showing that there are no polynomial-time algorithms that can always yield a solution that is within any constant factor of the optimal (unless P $=$ NP), i.e., for all polynomial-time algorithms $\mathcal{A}$ (resp., $\mathcal{A}'$) and $\forall K\in\mathbb{R}_{\ge 1}$, there are instances of the priori (resp., posteriori) KFSS problem where $r_{\mathcal{A}}(\Sigma)>K$ (resp., $r_{\mathcal{A}'}(\Sigma^*)>K$). 
\begin{remark}
Note that the ``constant'' in ``constant-factor approximation algorithm'' refers to the fact that the {\it cost} of the solution provided by the algorithm is upper-bounded by some (system-independent) constant times the cost of the optimal solution.  The algorithm can, however, use the system parameters when finding the solution.  For example, an optimal algorithm for the KFSS problem will be a $1$-factor approximation, and would use the system matrices, sensor costs, and budget to find the optimal solution.  Similarly, a polynomial-time $K$-factor approximation algorithm for KFSS would use the system parameters to produce a solution whose cost is guaranteed to be no more than $K$  times the cost of the optimal solution.  As indicated above, we will show that no such algorithm exists for any constant $K$ (unless P $=$ NP).
\end{remark}

To show the inapproximability of the priori KFSS problem, we relate it to the EXACT COVER BY 3-SETS $(X3C)$ problem described below \cite{garey1979computers}.

\begin{definition}
\label{def:X3C}
$(X3C)$ Given a finite set $D=\{d_1,\dots,d_{3m}\}$ and a collection $\mathcal{C}=\{c_1,\dots,c_{\tau}\}$ of $3$-element subsets of $D$, an {\it exact cover} for $D$ is a subcollection $\mathcal{C}'\subseteq \mathcal{C}$ such that every element of $D$ occurs in exactly one member of $\mathcal{C}'$.
\end{definition}

We will use the following result \cite{garey1979computers}.

\begin{lemma}
Given a finite set $D=\{d_1,\dots,d_{3m}\}$ and a collection $\mathcal{C}=\{c_1,\dots,c_{\tau}\}$ of $3$-element subsets of $D$, the problem of determining whether $\mathcal{C}$ contains an exact cover for $D$ is NP-complete.
\label{lemma:X3C}
\end{lemma}

As argued in Remark $\ref{remark:special instances}$, in order to show that the priori KFSS problem cannot be approximated within any constant factor in polynomial time, it is sufficient for us to show that certain special instances of this problem are inapproximable. Specifically, consider any instance of the $X3C$ problem. Using the results in Lemma $\ref{Lemma:minimum trace of sigma}$-$\ref{lemma:estimation error of state 1}$,  we will first construct an instance of the priori KFSS problem in polynomial time such that the difference between the solution to KFSS when the answer to $X3C$ is ``yes'' and the solution to KFSS when the answer to $X3C$ is ``no'' is large enough. Thus, we can then apply any hypothetical polynomial-time constant-factor approximation algorithm for the priori KFSS problem to the constructed priori KFSS instance and obtain the answer to the $X3C$ instance. Since we know from Lemma $\ref{lemma:X3C}$ that the $X3C$ problem is NP-complete, we obtain from Lemma $\ref{lemma:no polynomial algo}$ the following result; the detailed proof is provided in Appendix B.

\begin{theorem}
If P $\neq$ NP, then there is no  polynomial-time constant-factor approximation algorithm for the priori KFSS problem.
\label{thm:inapprox KFSS}
\end{theorem}

The following result is a direct consequence of the above arguments; the proof is also provided in Appendix B.
\begin{corollary}
\label{coro:inapprox KFSS post}
If P $\neq$ NP, then there is no polynomial-time constant-factor  approximation algorithm for the posteriori KFSS problem.
\end{corollary} 

\subsection{Inapproximability of the KFSA Problem}
In this section, we analyze the achievable performance of algorithms for the priori and posteriori KFSA problems. For any given algorithm $\mathcal{A}$ (resp., $\mathcal{A}'$) for the priori (resp., posteriori) KFSA problem, we define the following ratios:
\begin{equation}
r_{\mathcal{A}}(\tilde{\Sigma})\triangleq\frac{\text{trace}(\tilde{\Sigma}_{opt})}{\text{trace}(\tilde{\Sigma}_{\mathcal{A}})},
\label{eqn:performance ratio of KFSA algorithms}
\end{equation} 
and
\begin{equation}
r_{\mathcal{A}'}(\tilde{\Sigma}^*)\triangleq\frac{\text{trace}(\tilde{\Sigma}^*_{opt})}{\text{trace}(\tilde{\Sigma}^*_{\mathcal{A}'})},
\label{eqn:performance ratio of post KFSA algorithms}
\end{equation} 
where $\tilde{\Sigma}_{opt}$ (resp., $\tilde{\Sigma}^*_{opt}$) is the optimal solution to the priori (resp., posteriori) KFSA problem and $\tilde{\Sigma}_{\mathcal{A}}$ (resp., $\tilde{\Sigma}^*_{\mathcal{A}'}$) is the solution to the priori (resp., posteriori) KFSA problem given by algorithm $\mathcal{A}$ (resp., $\mathcal{A}'$). It is worth noting that using the arguments in \cite{zhang2017sensor}, the same (system-dependent) upper bounds for $r_{\mathcal{A}}(\tilde{\Sigma})$ and $r_{\mathcal{A}'}(\tilde{\Sigma}^*)$ can be obtained as those for $r_{\mathcal{A}}(\Sigma)$ and $r_{\mathcal{A}'}(\Sigma^*)$ in \cite{zhang2017sensor}, respectively. Nevertheless, we show that (if P$\neq$ NP) there is again no polynomial-time constant-factor  approximation algorithm for the priori (resp., posteriori) KFSA problem, i.e., for all  $K\in\mathbb{R}_{\ge1}$ and for all polynomial-time algorithms $\mathcal{A}$ (resp., $\mathcal{A}'$), there are instances of the priori (resp., posteriori) KFSA problem where $r_{\mathcal{A}}(\tilde{\Sigma})>K$ (resp., $r_{\mathcal{A'}}(\tilde{\Sigma}^*)>K$). To establish this result, we relate the priori KFSA problem to the $X3C$ problem described in Definition $\ref{def:X3C}$ and Lemma $\ref{lemma:X3C}$. Similarly to the proof of Theorem $\ref{thm:inapprox KFSS}$, given any instance of $X3C$, we will construct an instance of the priori KFSA problem and show that any hypothetical polynomial-time constant-factor approximation algorithm for the priori KFSA problem can be used to solve the $X3C$ problem. This leads to the following result; the detailed proof is provided in Appendix C. 
\begin{theorem}
\label{thm:inapprox KFSA}
If P $\neq$ NP, then there is no polynomial-time constant-factor approximation algorithm for the priori KFSA problem.
\end{theorem}

The arguments above also imply the following result whose proof is provided in Appendix C.
\begin{corollary}
\label{coro:inapprox KFSA post}
If P $\neq$ NP, then there is no  polynomial-time constant-factor approximation algorithm for the posteriori KFSA problem.
\end{corollary}

\section{Failure of Greedy Algorithms}\label{sec:greedy examples}
Our results in Theorem~\ref{thm:inapprox KFSS} and Theorem~\ref{thm:inapprox KFSA} indicate that no polynomial-time algorithm can be guaranteed to yield a solution that is within any constant factor of the optimal solution to the priori (resp., posteriori) KFSS and KFSA problems. In particular, these results apply to the greedy algorithms that are often studied for sensor selection in the literature (e.g., \cite{zhang2017sensor}, \cite{laszka2015resilient}), where sensors are iteratively selected (resp., attacked) in order to produce the greatest decrease (resp., increase) in the error covariance at each iteration. In this section we will focus on such greedy algorithms for the priori (resp., posteriori) KFSS and KFSA problems, and show explicitly how these greedy algorithms can fail to provide good solutions even for small and fixed instances with only three states; this provides additional insight into the factors that cause the KFSS and KFSA problems to be challenging.

\subsection{Failure of Greedy Algorithms for the KFSS Problem}
It was shown via simulations in \cite{zhang2017sensor} that greedy algorithms for KFSS work well in practice (e.g., for randomly generated systems).  In this section, we provide an explicit example showing that greedy algorithms for  the priori and posteriori KFSS problems can perform arbitrarily poorly, even for small systems (containing only three states). We consider the greedy algorithm for the priori (resp., posteriori) KFSS problem given in Algorithm \ref{algorithm:greedy KFSS}, for instances where all sensors have selection costs equal to $1$, and the sensor  selection budget $B\in\{1,\dots,q\}$ (i.e., up to $B$ sensors can be chosen).  For any such instance of the priori (resp., posteriori) KFSS problem, define $r_{gre}(\Sigma)=\frac{\textrm{trace}(\Sigma_{gre})}{\textrm{trace}(\Sigma_{opt})}$ (resp., $r_{gre}(\Sigma^*)=\frac{\textrm{trace}(\Sigma^*_{gre})}{\textrm{trace}(\Sigma^*_{opt})}$), where $\Sigma_{gre}$ (resp., $\Sigma^*_{gre}$) is the solution of Eq. \eqref{eqn:DARE} (resp., Eq. \eqref{eqn:postDARE}) corresponding to the sensors selected by Algorithm \ref{algorithm:greedy KFSS}.

\begin{algorithm}[H]
\textbf{Input:} An instance of priori (resp., posteriori) KFSS\\
\textbf{Output:} A set $\mathcal{S}$ of selected sensors
\caption{Greedy Algorithm for Problem $\ref{pro:priori KFSS problem}$}\label{algorithm:greedy KFSS}
\begin{algorithmic}[1]
\State $k\gets 1$, $\mathcal{S}\gets\emptyset$
\For{$k\le B$}
    \State $j\in\mathop{\arg\min}_{i\notin\mathcal{S}}\textrm{trace}(\Sigma(\mathcal{S}\cup\{i\}))$ (resp., $j\in\mathop{\arg\min}_{i\notin\mathcal{S}}\textrm{trace}(\Sigma^*(\mathcal{S}\cup\{i\}))$)
    \State $\mathcal{S}\gets \mathcal{S}\cup\{j\}$, $k\gets k+1$
\EndFor
\end{algorithmic}
\end{algorithm}

\begin{example}
Consider an instance of the priori (resp., posteriori) KFSS problem with matrices $W=I_3$ and $V=\mathbf{0}_{3\times 3}$, and $A$, $C$ defined as
\begin{equation*}
\label{eq:example matrices KFSS}
  A=
  \begin{bmatrix}
    \lambda_1 & 0 & 0 \\
            0 & 0 & 0\\
            0 & 0 & 0
  \end{bmatrix},\
 C=
  \begin{bmatrix}
    1 & h & h \\
    1 & 0 & h \\
    0 & 1 & 1
  \end{bmatrix},
\end{equation*}
where $0<|\lambda_1|<1$, $\lambda_1\in\mathbb{R}$, and $h\in\mathbb{R}_{>0}$. In addition, we have the selection budget $B=2$, the cost vector $b=[1\ 1\ 1]^T$ and the set of candidate sensors $\mathcal{Q}=\{1,2,3\}$, where sensor $i$ corresponds to the $i$th row of matrix $C$, for $i\in\{1,2,3\}$.
\label{exp:greedy bad results KFSS}
\end{example}

Based on the system defined in Example $\ref{exp:greedy bad results KFSS}$, we have the following result whose proof is provided in Appendix D.
\begin{theorem}
For the instance of the priori (resp., posteriori) KFSS problem defined in Example $\ref{exp:greedy bad results KFSS}$, the ratios $r_{gre}(\Sigma)=\frac{\textrm{trace}(\Sigma_{gre})}{\textrm{trace}(\Sigma_{opt})}$ and $r_{gre}(\Sigma^*)=\frac{\textrm{trace}(\Sigma^*_{gre})}{\textrm{trace}(\Sigma^*_{opt})}$ satisfy
\begin{equation}
\mathop{\lim}_{h\to\infty}r_{gre}(\Sigma)= \frac{2}{3}+\frac{1}{3(1-\lambda_1^2)},
\label{eqn:limit ratio KFSS}
\end{equation}
and
\begin{equation}
\mathop{\lim}_{h\to\infty}r_{gre}(\Sigma^*)= \frac{1}{1-\lambda_1^2},
\label{eqn:limit ratio post KFSS}
\end{equation}
respectively.
\label{thm:ratio of special system KFSS}
\end{theorem}

Examining Eq. \eqref{eqn:limit ratio KFSS} (resp., Eq. \eqref{eqn:limit ratio post KFSS}), we see that for the given instance of the priori (resp., posteriori) KFSS problem, we have $r_{gre}(\Sigma) \to \infty$ (resp., $r_{gre}(\Sigma^*)\to \infty$) as $h \to \infty$ and $\lambda_1 \to1$.  Thus, $r_{gre}(\Sigma)$ (resp., $r_{gre}(\Sigma^*)$) can be made arbitrarily large by choosing the parameters in the instance appropriately.  To explain the result in Theorem $\ref{thm:ratio of special system KFSS}$, we first note that the only nonzero eigenvalue of the diagonal $A$ defined in Example $\ref{exp:greedy bad results KFSS}$ is $\lambda_1$, and so we know from Lemma $\ref{Lemma:minimum trace of sigma}(a)$ that state $2$ and state $3$ of the system defined in Example $\ref{exp:greedy bad results KFSS}$ each contributes at most $1$ to $\text{trace}(\Sigma(\mu))$ (resp., $\text{trace}(\Sigma^*(\mu))$) for all $\mu$. Hence, in order to minimize $\text{trace}(\Sigma(\mu))$ (resp., $\text{trace}(\Sigma^*(\mu))$), we need to minimize the MSEE of state $1$. Moreover, the measurements of state $2$ and state $3$ can be viewed as measurement noise that corrupts the measurements of state $1$. It is then easy to observe from the form of matrix $C$ defined in Example $\ref{exp:greedy bad results KFSS}$ that sensor $2$ is the single best sensor among the three sensors since it provides measurements of state $1$ with less noise than sensor $1$ (and sensor $3$ does not measure state $1$ at all). Thus, the greedy algorithm for the priori (resp., posteriori) KFSS problem defined in Algorithm $\ref{algorithm:greedy KFSS}$ selects sensor $2$ in its first iteration. Nonetheless, we notice from $C$ defined in Example $\ref{exp:greedy bad results KFSS}$ that the optimal set of two sensors that minimizes $\text{trace}(\Sigma(\mu))$ (resp., $\text{trace}(\Sigma^*(\mu))$) contains sensor $1$ and sensor $3$, which together give us exact measurements (without measurement noise) on state $1$ (after some elementary row operations). Since the greedy algorithm selects sensor $2$ in its first iteration, no matter which sensor it selects in its second iteration, the two chosen sensors can only give a noisy measurement of state $1$ (if we view the measurements of state $2$ and state $3$ as measurement noise), and the variance of the measurement noise can be made arbitrary large if we take $h\to\infty$ in $C$ defined in Example $\ref{exp:greedy bad results KFSS}$.  Hence, the greedy algorithm fails to perform well due to its myopic choice in the first iteration.

It is also useful to note that the above behavior holds for any algorithm that outputs a sensor selection that contains sensor $2$ for the above example.

\subsection{Failure of Greedy Algorithms for the KFSA Problem}
In \cite{laszka2015resilient}, the authors showed that a simple greedy algorithm can perform arbitrarily poorly for an instance of the observation attack problem in Gaussian process regression. Here, we consider a simple greedy algorithm for the priori (resp., posteriori) KFSA problem given in Algorithm \ref{algorithm: greedy KFSA}, for instances where all sensors have an attack cost of $1$, and the sensor attack budget $\Omega\in\{1,\dots,q\}$ (i.e., up to $\Omega$ sensors can be attacked).  For any such instance of the priori (resp., posteriori) KFSA problem, define $r_{gre}(\tilde{\Sigma})=\frac{\textrm{trace}(\tilde{\Sigma}_{opt})}{\textrm{trace}(\tilde{\Sigma}_{gre})}$ (resp., $r_{gre}(\tilde{\Sigma}^*)=\frac{\textrm{trace}(\tilde{\Sigma}^*_{opt})}{\textrm{trace}(\tilde{\Sigma}^*_{gre})}$), where $\tilde{\Sigma}_{gre}$ (resp., $\tilde{\Sigma}^*_{gre}$) is the solution to the priori (resp., posteriori) KFSA problem given by Algorithm \ref{algorithm: greedy KFSA}. We then show that Algorithm \ref{algorithm: greedy KFSA} can perform arbitrarily poorly for a simple instance of the priori (resp., posteriori) KFSA problem  described below.

\begin{algorithm}
\textbf{Input:} An instance of priori (resp., posteriori) KFSA\\
\textbf{Output:} A set $\mathcal{S}$ of targeted sensors
\caption{Greedy Algorithm for Problem $\ref{pro:priori KFSA problem}$}\label{algorithm: greedy KFSA}
\begin{algorithmic}[1]
\State $k\gets 1$, $\mathcal{S}\gets\emptyset$
\For{$k\le \Omega$}
    \State $j\in\mathop{\arg\max}_{i\notin\mathcal{S}}\textrm{trace}(\Sigma(\mathcal{Q}\setminus(\mathcal{S}\cup\{i\})))$ (resp., $j\in\mathop{\arg\max}_{i\notin\mathcal{S}}\textrm{trace}(\Sigma^*(\mathcal{Q}\setminus(\mathcal{S}\cup\{i\})))$)
    \State $\mathcal{S}\gets \mathcal{S}\cup\{j\}$, $k\gets k+1$
\EndFor
\end{algorithmic}
\end{algorithm}

\begin{example}
Consider an instance of the priori (resp., posteriori) KFSA problem with matrices $W=I_3$, $V=\mathbf{0}_{4\times 4}$, and $A$, $C$ defined as
\begin{equation*}
\label{eq:example matrices KFSA}
  A=
  \begin{bmatrix}
    \lambda_1 & 0 & 0 \\
            0 & 0 & 0\\
            0 & 0 & 0
  \end{bmatrix},\
  C=
  \begin{bmatrix}
    1 & h& h \\
    1 & 0 & h\\
    0 & 1 & 0 \\
    0 & 0 & 1
  \end{bmatrix},\\
\end{equation*}
where $0<|\lambda_1|<1$, $\lambda_1\in\mathbb{R}$ and $h\in\mathbb{R}_{>0}$. In addition, the attack budget is $\Omega=2$, the cost vector is $\omega=[1\ 1\ 1\ 1]^T$, and the set of sensors $\mathcal{Q}=\{1,2,3,4\}$ has already been installed on the system, where sensor $i$ corresponds to the $i$th row of matrix $C$, for $i\in\{1,2,3,4\}$.
\label{exp:greedy bad results KFSA}
\end{example}

We then have the following result, whose proof is provided in Appendix D.
\begin{theorem}
For the instance of the priori (resp., posteriori) KFSA problem defined in Example \ref{exp:greedy bad results KFSA}, the ratios $r_{gre}(\tilde{\Sigma})=\frac{\textrm{trace}(\tilde{\Sigma}_{opt})}{\textrm{trace}(\tilde{\Sigma}_{gre})}$ and $r_{gre}(\tilde{\Sigma}^*)=\frac{\textrm{trace}(\tilde{\Sigma}^*_{opt})}{\textrm{trace}(\tilde{\Sigma}^*_{gre})}$ satisfy
\begin{equation}
\mathop{\lim}_{h\to 0}r_{gre}(\tilde{\Sigma})= \frac{2}{3}+\frac{1}{3(1-\lambda_1^2)},
\label{eqn:limit ratio KFSA}
\end{equation}
and
\begin{equation}
\mathop{\lim}_{h\to 0}r_{gre}(\tilde{\Sigma}^*)= \frac{1}{1-\lambda_1^2},
\label{eqn:limit ratio post KFSA}
\end{equation}
\label{thm:ratio of special system KFSA}
respectively.
\end{theorem}

Inspecting Eq. \eqref{eqn:limit ratio KFSA} (resp., Eq. \eqref{eqn:limit ratio post KFSA}), we observe that for the given instance of the priori (resp., posteriori) KFSA problem, we have $r_{gre}(\tilde{\Sigma}) \to \infty$ (resp., $r_{gre}(\tilde{\Sigma}^*)\to \infty$) as $h \to 0$ and $\lambda_1 \to1$.  Thus, $r_{gre}(\tilde{\Sigma})$ (resp., $r_{gre}(\tilde{\Sigma}^*)$) can be made arbitrarily large by choosing the parameters in the instance appropriately. Here, we explain the results in Theorem $\ref{thm:ratio of special system KFSA}$ as follows. Using similar arguments to those before, we know from the structure of matrix $A$ defined in Example $\ref{exp:greedy bad results KFSA}$ that in order to maximize $\text{trace}(\Sigma(\nu^c))$ (resp., $\text{trace}(\Sigma^*(\nu^c))$), we need to maximize the MSEE of state $1$, i.e., make the measurements of state $1$ ``worse''. Again, the measurements of state $2$ and state $3$ can be viewed as measurement noise that corrupts the measurements of state $1$. No matter which of sensor $1$, sensor $2$, or sensor $3$ is attacked, the resulting measurement matrix $C(\nu^c)$ is full column rank, which yields an exact measurement of state $1$. We also observe that if sensor $4$ is targeted, the surviving sensors can only provide measurements of state $1$ that are corrupted by measurements of states $2$ and state $3$. Hence, the greedy algorithm for the priori (resp., posteriori) KFSA problem defined in Algorithm $\ref{algorithm: greedy KFSA}$ targets sensor $4$ in its first iteration, since it is the single best sensor to attack from the four sensors. Nevertheless, sensor $1$ and sensor $2$ form the optimal set of sensors to be attacked to maximize $\text{trace}(\Sigma(\nu^c))$ (resp., $\text{trace}(\Sigma^*(\nu^c))$), since the surviving sensors provide no measurement of state $1$. Since the greedy algorithm targets sensor $4$ in its first iteration, no matter which sensor it targets in the second step, the surviving sensors can always provide some measurements of state $1$ with noise (if we view the measurements of state $2$ and state $3$ as measurement noise), and the variance of the noise will vanish if we take $h\to 0$ in matrix $C$ defined in Example $\ref{exp:greedy bad results KFSA}$. Hence, the myopic behavior of the greedy algorithm makes it perform poorly.

Furthermore, it is useful to note that the above result holds for any algorithm that outputs a sensor attack that does not contain sensor $1$ or sensor $2$ for the above example.

\begin{remark}
Using similar arguments to those in the proof of Theorem $\ref{thm:ratio of special system KFSS}$ (resp., Theorem $\ref{thm:ratio of special system KFSA}$), we can also show that when we set $V=\varepsilon I_{3}$ (resp., $V=\varepsilon I_{4}$), where $\varepsilon\in\mathbb{R}_{>0}$, the results in Eqs. \eqref{eqn:limit ratio KFSS}-\eqref{eqn:limit ratio post KFSS} (resp., Eqs. \eqref{eqn:limit ratio KFSA}-\eqref{eqn:limit ratio post KFSA}) hold if we let $\varepsilon\to0$. This phenomenon  is also observed in \cite{chamon2020approximate}, where the approximation guarantees for the greedy algorithms provided in that paper get worse as the sensor measurement noise tends to zero.
\end{remark}

\section{Conclusions}\label{sec:conclusion}

In this paper, we studied sensor selection and attack problems for (steady state) Kalman filtering of linear dynamical systems. We showed that these problems are NP-hard and have no polynomial-time constant-factor  approximation algorithms, even under the assumption that the system is stable and each sensor has identical cost. To illustrate this point, we provided explicit examples showing how greedy algorithms can perform arbitrarily poorly on these problems, even when the system only has three states. Our results shed new insights into the problem of sensor selection and attack for Kalman filtering and show, in particular, that this problem is more difficult than other variants of the sensor selection problem that have submodular (or supermodular) cost functions. Future work on extending the results to Kalman filtering over finite time horizons, characterizing achievable (non-constant) approximation ratios, identifying classes of systems that admit near-optimal approximation algorithms, and investigating resilient sensor selection problems under adversarial settings would be of interest.

\section{Acknowledgments}
The authors thank the anonymous reviewers for their insightful comments that helped to improve the paper.


%%%%%%%%%%%%%%%%%%%%%%%%%%%%%%%%%%%%%%%%%%%%%%%%%%%%%%%%%%%%%%%%%%%%%%%%%%%%%%%%

%%%%%%%%%%%%%%%%%%%%%%%%%%%%%%%%%%%%%%%%%%%%%%%%%%%%%%%%%%%%%%%%%%%%%%%%%%%%%%%%

%%%%%%%%%%%%%%%%%%%%%%%%%%%%%%%%%%%%%%%%%%%%%%%%%%%%%%%%%%%%%%%%%%%%%%%%%%%%%%%%
\section*{Appendix A}

\subsection*{Proof of Lemma \ref{Lemma:minimum trace of sigma}:}
Since $A$ and $W$ are diagonal, the system represents a set of $n$ scalar subsystems of the form 
\begin{equation*}
x_i[k+1]=\lambda_i x_i[k] + w_i[k], \forall i\in\{1,\dots,n\},
\label{eqn:decoupled system dynamics}
\end{equation*}
where $x_i[k]$ is the $i$th state of $x[k]$ and $w_i[k]$ is a zero-mean white Gaussian noise process with variance $\sigma_{w_i}^2=W_{ii}$. As $A$ is stable, the pair $(A,C(\mu))$ is detectable and the pair $(A,W^{\frac{1}{2}})$ is stabilizable for all sensor selections $\mu$. Thus, the limits $\displaystyle\mathop{\lim}_{k\to\infty}(\Sigma_{k|k-1}(\mu))_{ii}$ and $\displaystyle\mathop{\lim}_{k\to\infty}(\Sigma_{k|k}(\mu))_{ii}$ exist for all $i$ and for all $\mu$ (based on Lemma \ref{lemma:Anderson optimal filtering}), and are denoted as $(\Sigma(\mu))_{ii}$ and $(\Sigma^*(\mu))_{ii}$, respectively.

Proof of (a): Since $A$ and $W$ are diagonal, we know from Eq. $(\ref{eqn:couple priori and posteriori 1})$ that 
\begin{equation*}
(\Sigma(\mu))_{ii}=\lambda_i^2(\Sigma^*(\mu))_{ii}+W_{ii},
\label{eqn:coupled priori and posteriori diagonal case}
\end{equation*}
which implies $(\Sigma(\mu))_{ii}\ge W_{ii}$, $\forall i\in\{1,\dots,n\}$. Moreover, it is easy to see that $(\Sigma(\mu))_{ii}\le (\Sigma(\mathbf{0}))_{ii}$, $\forall i\in\{1,\dots,n\}$. Since $C(\mathbf{0})=\mathbf{0}$, we obtain from Eq. $(\ref{eqn:DARE})$ that
\begin{equation*}
\Sigma(\mathbf{0})=A\Sigma(\mathbf{0})A^T+W.
\end{equation*}
which implies that $(\Sigma(\mathbf{0}))_{ii}=\frac{W_{ii}}{1-\lambda_i^2}$ since $A$ is diagonal. Hence, $W_{ii}\le(\Sigma(\mu))_{ii}\le\frac{W_{ii}}{1-\lambda_i^2}$, $\forall i\in\{1,\dots,n\}$. Similarly, we also have $(\Sigma^*(\mu))_{ii}\le(\Sigma^*(\mathbf{0}))_{ii}$, and we obtain from Eq. \eqref{eqn:postDARE} that
\begin{equation*}
\Sigma^*(\mathbf{0})=A\Sigma^*(\mathbf{0})A^T+W.
\end{equation*}
Thus, $0\le(\Sigma^*(\mu))_{ii}\le\frac{W_{ii}}{1-\lambda_i^2}$, $\forall i\in\{1,\dots,n\}$.

Proof of (b): Assume without loss of generality that the first column of $C(\mu)$ is zero, since we can simply renumber the states to make this the case without affecting the trace of the error covariance matrix. We then have  $C(\mu)$ of the form
\begin{equation*}
C(\mu)=\begin{bmatrix}
\mathbf{0} & C_1(\mu)
\end{bmatrix}.
 \end{equation*}
Moreover, since $A$ and $W$ are diagonal, we obtain from Eq. $(\ref{eqn:DARE})$ that $\Sigma(\mu)$ is of the form
\begin{equation*}
\Sigma(\mu)=\begin{bmatrix}
\Sigma_1(\mu) & \mathbf{0}\\
\mathbf{0} & \Sigma_2(\mu)
\end{bmatrix},
\end{equation*}
where $\Sigma_1(\mu)=(\Sigma(\mu))_{11}$ and satisfies 
\begin{equation*}
(\Sigma(\mu))_{11}=\lambda_i^2(\Sigma(\mu))_{11}+W_{11}, 
\end{equation*}
which implies $(\Sigma(\mu))_{11}=\frac{W_{11}}{1-\lambda_1^2}$. Furthermore, we obtain from Eq. \eqref{eqn:postDARE} that $\Sigma^*(\mu)$ is of the form
\begin{equation*}
\Sigma^*(\mu)=\begin{bmatrix}
\Sigma^*_1(\mu) & \mathbf{0}\\
\mathbf{0} & \Sigma^*_2(\mu)
\end{bmatrix},
\end{equation*}
where $\Sigma^*_1(\mu)=(\Sigma^*(\mu))_{11}$ and satisfies 
\begin{equation*}
(\Sigma^*(\mu))_{11}=\lambda_1^2(\Sigma^*(\mu))_{11}+W_{11}, 
\end{equation*}
which implies $(\Sigma^*(\mu))_{11}=\frac{W_{11}}{1-\lambda_1^2}$.

Proof of (c): We assume without loss of generality that $\mathbf{e}_{1}\in\textrm{rowspace}(C(\mu))$. If we further perform elementary row operations on $C(\mu)$, which does not change the solution to Eq. \eqref{eqn:DARE} (resp., Eq. \eqref{eqn:postDARE}), we obtain a measurement matrix $\tilde{C}(\mu)$ of the form 
\begin{equation*}
\tilde{C}(\mu)=\begin{bmatrix}
1 & \mathbf{0}\\
\mathbf{0} & \tilde{C}_1(\mu)
\end{bmatrix}
 \end{equation*}
 with $\tilde{V}(\mu)=\mathbf{0}$. Moreover, since $A$ and $W$ are diagonal, we obtain from Eq. $(\ref{eqn:DARE})$ that $\Sigma(\mu)$ is of the form
\begin{equation*}
\Sigma(\mu)=\begin{bmatrix}
\Sigma_1(\mu) & \mathbf{0}\\
\mathbf{0} & \Sigma_2(\mu)
\end{bmatrix},
\end{equation*}
where $\Sigma_1(\mu)=(\Sigma(\mu))_{11}$ and satisfies
\begin{equation*}
(\Sigma(\mu))_{11}=\lambda_1^2(\Sigma(\mu))_{11}+W_{11}-\lambda_1^2(\Sigma(\mu))_{11},
\end{equation*}
which implies $(\Sigma(\mu))_{11}=W_{11}$. Furthermore, we obtain from Eq. $(\ref{eqn:postDARE})$ that $\Sigma^*(\mu)$ is of the form 
\begin{equation*}
\Sigma^*(\mu)=\begin{bmatrix}
\Sigma^*_1(\mu) & \mathbf{0}\\
\mathbf{0} & \Sigma^*_2(\mu)
\end{bmatrix},
\end{equation*}
where $\Sigma^*_1(\mu)=(\Sigma^*(\mu))_{11}$ and satisfies $(\Sigma^*(\mu))_{11}=(\Sigma(\mu))_{11}-(\Sigma(\mu))_{11}=0$.\hfill\IEEEQED

\subsection*{Proof of Lemma \ref{lemma:estimation error of state 1}:}
Proof of (a): We first note from Lemma $\ref{lemma:Anderson optimal filtering}$ that the limit $\Sigma(\mu)$ exists for all $\mu$ (since $A$ is stable). Since $A=\textrm{diag}(\lambda_1,0,\dots,0)$, we have $x_i[k+1]=w_i[k]$, $\forall i\in\{2,\dots,n\}$ and $\forall k\in\mathbb{Z}_{\ge0}$. Moreover, we have from Eq. $\eqref{eqn:all sensors measurements}$ that
\begin{equation*}
y[k]=[1\ \mathbf{0}_{1\times (n-1)}]x[k]+v[k]+v^{\prime}[k]=x_1[k]+\tilde{v}[k],\ \forall k\in\mathbb{Z}_{\ge0},
\end{equation*}
where $v^{\prime}[k]=\displaystyle\sum_{i=1}^{n-1}\gamma_i x_{i+1}[k]$ and $\tilde{v}[k]\triangleq v[k]+v^{\prime}[k]$. Recall that we have assumed with out loss of generality that $\bar{x}_0=\mathbf{0}$ and $\Pi_0=I_n$. Moreover, noting that $W=I_n$ and that $x[0]$ is independent of $w[k]$ and $v[k]$ for all $k\in\mathbb{Z}_{\ge0}$, where $w[k]$ and $v[k]$ are uncorrelated zero-mean white Gaussian noise processes (as assumed), we have that $\tilde{v}[k]$ is a zero-mean white Gaussian noise process with $\mathbb{E}[(\tilde{v}[k])^2]=\lVert\gamma\rVert_2^2+\sigma_v^2$. Thus, to compute the MSEE of state $1$ of the Kalman filter, i.e., $\Sigma_{11}$, we can consider a scalar discrete-time linear system with $A=\lambda_1$, $C=1$, $W=1$ and $V=\alpha^2$, and obtain from Eq. $\eqref{eqn:DARE}$ the scalar DARE 
\begin{equation}
\Sigma_{11}=\lambda_1^2(1-\frac{\Sigma_{11}}{\alpha^2+\Sigma_{11}})\Sigma_{11}+1,
\label{eqn:scalar DARE for state 1}
\end{equation}
where $\alpha^2=\|\gamma\|^2_2+\sigma_v^2$. Solving for $\Sigma_{11}$ in Eq. $(\ref{eqn:scalar DARE for state 1})$ and omitting the negative solution lead to Eq. $(\ref{eqn: mean square estimation error of state 1})$. 

To show that $\Sigma_{11}$ is strictly increasing in $\alpha^2\in\mathbb{R}_{\ge0}$, we can use the result of Lemma $6$ in \cite{zhang2017sensor}. For a discrete-time linear system defined in Eq. $(\ref{eqn:system dynamics})$ and Eq. $(\ref{eqn:all sensors measurements})$, given $A=\lambda_1$ and $W=1$, suppose we have two sensors with the measurement matrices as $C_1=C_2=1$ and the variances of the (Gaussian) measurement noise as $V_1=\alpha_1^2$ and $V_2=\alpha_2^2$. Define $R\triangleq C^TV^{-1}C$ to be the sensor information matrix corresponding to a sensor with measurement matrix $C$ and measurement noise covariance matrix $V$. The sensor information matrix of these two sensors are denoted as $R_1$ and $R_2$. We then have $R_1=\frac{1}{\alpha_1^2}$ and $R_2=\frac{1}{\alpha_2^2}$. If $\alpha_1^2>\alpha_2^2$, we know from Lemma $6$ in \cite{zhang2017sensor} that $\Sigma_{11}(\alpha_1^2)<\Sigma_{11}(\alpha_2^2)$. Hence, $\Sigma_{11}(\alpha^2)$ is a strictly increasing function of $\alpha^2\in\mathbb{R}_{\ge0}$. For $\alpha>0$, we can rewrite Eq. $(\ref{eqn: mean square estimation error of state 1})$ as
\begin{equation}
\Sigma_{11}(\alpha^2)=\frac{2}{\sqrt[]{(1-\lambda_1^2-\frac{1}{\alpha^2})^2+\frac{4}{\alpha^2}}+1-\lambda_1^2-\frac{1}{\alpha^2}}.
\label{eqn:rewrite mean square estimation error of state 1}
\end{equation}
We then obtain from Eq. $(\ref{eqn:rewrite mean square estimation error of state 1})$ that $\displaystyle\mathop{\lim}_{\alpha\to\infty}\Sigma_{11}(\alpha^2)=\frac{1}{1-\lambda_1^2}$.

Proof of (b): Using similar arguments to those above, we obtain from Eq. $\eqref{eqn:all sensors measurements}$ that
\begin{equation*}
y[k]=\mathbf{1}_{n-1}x_1[k]+v'[k],
\end{equation*}
where $v'[k]=\rho\left[\begin{matrix}x_2[k] \cdots x_n[k]\end{matrix}\right]^T$ is a zero-mean white Gaussian noise process with $\mathbb{E}[v'[k](v'[k])^T]=\rho^2 I_{n-1}$. Hence, to compute the MSEE of state $1$ of the Kalman filter, i.e., $\Sigma'_{11}$, we can consider a system with $A=\lambda_1$, $C=\mathbf{1}_{n-1}$, $W=1$ and $V=\rho^2 I_{n-1}$. Solving Eq. $\eqref{eqn:DARE}$ (using the matrix inversion lemma \cite{horn1985matrix}) yields Eq. $\eqref{eqn:mean square estimation error of state 1 with identity noise states}$. Similarly, we have $\Sigma'_{11}$ is strictly increasing in $\rho^2\in\mathbb{R}_{\ge0}$ and $\displaystyle\mathop{\lim}_{\rho\to\infty}\Sigma'_{11}(\rho^2)=\frac{1}{1-\lambda_1^2}$.\hfill\IEEEQED

\section*{Appendix B}
We will use the following result in the proof of Theorem $\ref{thm:inapprox KFSS}$.
\begin{lemma}
Consider an instance of $X3C$: a finite set $D$ with $|D|=3m$, and a collection $\mathcal{C}=\{c_1,\dots,c_{\tau}\}$ of $\tau$ $3$-element subsets of $D$, where $\tau\ge m$. For each element $c_i \in \mathcal{C}$, define a column vector $g_i \in \mathbb{R}^{3m}$ to encode which elements of $D$ are contained in $c_i$, i.e., for $i \in \{1, 2, \ldots, \tau\}$ and $j \in \{1, 2, \ldots, 3m\}$, $(g_i)_j = 1$ if element $j$ of set $D$ is in $c_i$, and $(g_i)_j = 0$ otherwise. Denote
$G\triangleq\left[\begin{matrix} g_1 & \cdots & g_{\tau}\end{matrix}\right]^T$. For any $l\le m$ ($l\in\mathbb{Z}$) and $\mathcal{L}\triangleq\{i_1,\dots,i_l\}\subseteq\{1,\dots,\tau\}$, define $G_{\mathcal{L}}\triangleq\left[\begin{matrix} g_{i_1} & \cdots & g_{i_l}\end{matrix}\right]^T$ and denote $\textrm{rank}(G_{\mathcal{L}})=r_{\mathcal{L}}$.\footnote{We drop the subscript $\mathcal{L}$ on $r$ for notational simplicity.}. If the answer to the $X3C$ problem is ``no'', then for all $\mathcal{L}$ with $|\mathcal{L}|\le m$, there exists an orthogonal matrix $N\in\mathbb{R}^{3m\times 3m}$ such that
\begin{equation}
\begin{bmatrix}
\mathbf{1}_{3m}^T \\
G_{\mathcal{L}}
\end{bmatrix}N=\begin{bmatrix}
\mathbf{\gamma} & \mathbf{\beta}\\
\mathbf{0} & \tilde{G}_{\mathcal{L}}
\end{bmatrix}, 
\label{eqn: X3C no solution}
\end{equation}
where $\tilde{G}_{\mathcal{L}}\in\mathbb{R}^{l\times r}$ is of full column rank, $\gamma\in\mathbb{R}^{1\times(3m-r)}$ and at least $\kappa\ge1$ ($\kappa\in\mathbb{Z}$) elements of $\mathbf{\gamma}$ are $1$'s , and $\beta\in\mathbb{R}^{1\times r}$. Further elementary row operations on $\left[\begin{smallmatrix}
\mathbf{\gamma} & \mathbf{\beta}\\
\mathbf{0} & \tilde{G}_{\mathcal{L}}
\end{smallmatrix}\right]$ transform it into the form $\left[\begin{smallmatrix}
\mathbf{\gamma} & \mathbf{0}\\
\mathbf{0} & \tilde{G}_{\mathcal{L}}
\end{smallmatrix}\right]$. 
\label{lemma:X3C no solution}
\end{lemma}
\begin{IEEEproof}
Assume without loss of generality that there are no identical subsets in $\mathcal{C}$. Since $\textrm{rank}(G_{\mathcal{L}})=r$, the dimension of the nullspace of $G_{\mathcal{L}}$ is $3m-r$. We choose an orthonormal basis of the nullspace of $G_{\mathcal{L}}$ and let it form the first $3m-r$ columns of $N$, denoted as $N_1$. Then, we choose an orthonormal basis of the columnspace of $G_{\mathcal{L}}^T$ and let it form the rest of the $r$ columns of $N$, denoted as $N_2$. Clearly, $N=\begin{bmatrix}N_1 & N_2\end{bmatrix}\in\mathbb{R}^{3m\times 3m}$ is an orthogonal matrix. Furthermore, since the answer to the $X3C$ problem is ``no'', for any union of $l\le m$ ($l\in\mathbb{Z}$) subsets in $\mathcal{C}$, denoted as $\mathcal{C}_l$, there exist $\kappa\ge1$ $(\kappa\in\mathbb{Z})$ elements in $D$ that are not covered by $\mathcal{C}_l$, i.e., $G_{\mathcal{L}}$ has $\kappa$ zero columns. Let these denote the $j_1$th, $\dots$, $j_{\kappa}$th columns of $G_{\mathcal{L}}$, where $\{j_1,\dots,j_{\kappa}\}\subseteq\{1,\dots,3m\}$. Hence, we can always choose $\mathbf{e}_{j_1},\dots,\mathbf{e}_{j_{\kappa}}$ to be in the orthonormal basis of the nullspace of $G_{\mathcal{L}}$, i.e., as columns of $N_1$. Constructing $N$ in this way, we have $G_{\mathcal{L}} N_1=\mathbf{0}$ and $G_{\mathcal{L}} N_2=\tilde{G}_{\mathcal{L}}$, where $\tilde{G}_{\mathcal{L}}\in\mathbb{R}^{l\times r}$ is of full column rank since the columns of $N_2$ form an orthonormal basis of the columnspace of $G^T_{\mathcal{L}}$ and $r\le l$. Moreover, we have $\mathbf{1}_{3m}^T N_1=\mathbf{\gamma}$ and $\mathbf{1}_{3m}^T N_2=\mathbf{\beta}$, where at least $\kappa$ elements of $\mathbf{\gamma}$ are $1$'s (since $\mathbf{1}_{3m}^T\mathbf{e}_{j_s}^T=1$, $\forall s\in\{1,\dots,\kappa\}$). Combining these results, we obtain Eq. $(\ref{eqn: X3C no solution})$. Since $\tilde{G}_{\mathcal{L}}$ is of full column rank, we can perform elementary row operations on $\left[\begin{smallmatrix}
\mathbf{\gamma} & \mathbf{\beta}\\
\mathbf{0} & \tilde{G}_{\mathcal{L}}
\end{smallmatrix}\right]$ and obtain $\left[\begin{smallmatrix}
\mathbf{\gamma} & \mathbf{0}\\
\mathbf{0} & \tilde{G}_{\mathcal{L}}
\end{smallmatrix}\right]$.
\end{IEEEproof}
\subsection*{Proof of Theorem $\ref{thm:inapprox KFSS}$:}
Assume that there exists such an approximation algorithm $\mathcal{A}$, i.e., $\exists K\in\mathbb{R}_{\ge 1}$ such that $r_{\mathcal{A}}(\Sigma)\le K$ for all instances of the priori KFSS problem, where $r_{\mathcal{A}}(\Sigma)$ is defined in Eq. $\eqref{eqn:performance ratio of algorithm}$. We will show that $\mathcal{A}$ can be used to solve the $X3C$ problem, which will lead to a contradiction. 

Given an arbitrary instance of the $X3C$ problem described in Definition $\ref{def:X3C}$ and Lemma $\ref{lemma:X3C}$, for each element $c_i \in \mathcal{C}$, we define $g_i \in \mathbb{R}^{3m}$ to encode which elements of $D$ are contained in $c_i$. Specifically, for $i \in \{1, 2, \ldots, \tau\}$ and $j \in \{1, 2, \ldots, 3m\}$, $(g_i)_j = 1$ if $d_j\in D$ is in $c_i$, and $(g_i)_j = 0$ otherwise. Denote $G\triangleq\left[\begin{matrix}g_1 & \cdots & g_{\tau}\end{matrix}\right]^T$. Thus $G^Tx=\mathbf{1}_{3m}$ has a solution $x\in\{0,1\}^{\tau}$ such that $x$ has $m$ nonzero entries if and only if the answer to the $X3C$ instance is ``yes'' \cite{natarajan1995sparse}. 

Given the above instance of $X3C$, we then construct an instance of the priori KFSS problem as follows. Denote $Z= \lceil K\rceil(m+1)(\sigma_v^2+3)$, where we set $\sigma_v=1$. Define the system dynamics matrix as $A=\textrm{diag}(\lambda_1,0,\dots,0)\in\mathbb{R}^{(3m+1)\times(3m+1)}$, where $\lambda_1=\frac{Z-1/2}{Z}$. Note that $Z\in\mathbb{Z}_{>1}$ and $0<\lambda_1<1$. The set $\mathcal{Q}$ is defined to contain $\tau+1$ sensors with collective measurement matrix
\begin{equation}
C=\begin{bmatrix}
1 & \varepsilon \mathbf{1}_{3m}^T\\
\mathbf{0} & G
\end{bmatrix},
\label{eqn: C matrix inapprox}
\end{equation}
where $G$ is defined based on the given instance of $X3C$ as above. The constant $\varepsilon$ is chosen as $\varepsilon=2Z\Bigl\lceil\sqrt{Z-1}\Bigr\rceil+1$. The system noise covariance matrix is set to be $W=I_{3m+1}$. The measurement noise covariance matrix is set as $V=\sigma_v^2\begin{bmatrix}1 & \mathbf{0}\\ \mathbf{0} & \frac{1}{\varepsilon^2}I_{\tau}\end{bmatrix}$. The sensor selection cost vector is set as $b=\mathbf{1}_{\tau+1}$, and the sensor selection budget is set as $B=m+1$.  Note that the sensor selection vector for this instance is denoted by $\mu\in\{0,1\}^{\tau+1}$. For the above construction, since the only nonzero eigenvalue of $A$ is $\lambda_1$, we know from Lemma $\ref{Lemma:minimum trace of sigma}$(a) that $\sum_{i=2}^{3m+1}(\Sigma(\mu))_{ii}=\sum_{i=2}^{3m+1}W_{ii}=3m$ for all $\mu$.

We claim that algorithm $\mathcal{A}$ will return a sensor selection vector $\mu$ such that $\textrm{trace}(\Sigma(\mu))\le K(m+1)(\sigma_v^2+3)$ if and only if the answer to the $X3C$ problem is ``yes''. 

We prove the above claim as follows. Suppose that the answer to the instance of the $X3C$ problem is ``yes''. Then $G^Tx=\mathbf{1}_{3m}$ has a solution such that $x$ has $m$ nonzero entries. Denote the solution as $x^*$ and denote $\textrm{supp}(x^*)=\{i_1,\dots,i_m\}$. Define $\tilde{\mu}$ to be the sensor selection vector that indicates selecting the first and the $(i_1+1)$th to the $(i_m+1)$th sensors, i.e., sensors that correspond to rows $C_1$, $C_{i_1+1},\dots, C_{i_m+1}$ from \eqref{eqn: C matrix inapprox}. Since $G^Tx^*=\mathbf{1}_{3m}$, we have $[1\ -\varepsilon x^{*T}]C=\mathbf{e}_1$ for $C$ defined in Eq. $(\ref{eqn: C matrix inapprox})$. Noting that $\textrm{supp}(x^*)=\{i_1,\dots,i_m\}$, it then follows that $\mathbf{e}_1\in\textrm{rowspace}(C(\tilde{\mu}))$. We can then perform elementary row operations on $C(\tilde{\mu})$ (which does not change the steady state {\it a priori} error covariance matrix of the corresponding Kalman filter) and obtain $\Gamma C(\tilde{\mu})\triangleq\tilde{C}(\tilde{\mu})=\begin{bmatrix}1 & \mathbf{0}\\ \mathbf{0} & *\end{bmatrix}$ with the corresponding measurement noise covariance $\Gamma V(\mu)\Gamma^T\triangleq\tilde{V}(\tilde{\mu})=\begin{bmatrix}\sigma_v^2(m+1) & *\\ * & *\end{bmatrix}$, where $\Gamma=\begin{bmatrix}1 & -\varepsilon\mathbf{1}^T_m \\ \mathbf{0} & I_m\end{bmatrix}$. Let $\tilde{\Sigma}$ denote the error covariance obtained from sensing matrix $(\tilde{C}(\tilde{\mu}))_1=\mathbf{e}_1$ with measurement noise variance $\tilde{\sigma}_v^2\triangleq\sigma_v^2(m+1)$, which corresponds to the first sensor in $\tilde{C}(\tilde{\mu})$. We then know from Lemma $\ref{lemma:estimation error of state 1}$(a) that
\begin{equation*}
\tilde{\Sigma}_{11}=\frac{1+\tilde{\sigma}_v^2\lambda_1^2-\tilde{\sigma}_v^2+\sqrt{(\tilde{\sigma}_v^2-\tilde{\sigma}_v^2\lambda_1^2-1)^2+4\tilde{\sigma}_v^2}}{2},
\end{equation*}
which further implies 
\begin{align}\nonumber
\tilde{\Sigma}_{11}&\le\frac{1+\sqrt{(\tilde{\sigma}_v^2(1-\lambda_1^2))^2-2\tilde{\sigma}_v^2(1-\lambda_1^2)+1+4\tilde{\sigma}_v^2}}{2}\\\nonumber
&\le\frac{1+\sqrt{(\tilde{\sigma}_v^2(1-\lambda_1^2))^2+1+4\tilde{\sigma}_v^2}}{2}\\
&\le\frac{1+\sqrt{\tilde{\sigma}_v^4+4\tilde{\sigma}_v^2+4}}{2}\le\frac{1+\tilde{\sigma}_v^2+2}{2}.
\label{eqn:thm1 upper bound on MSEE of state 1}
\end{align}
Using similar arguments to those above, we have that $\sum_{i=2}^{3m+1}\tilde{\Sigma}_{ii}=3m$. We then obtain from \eqref{eqn:thm1 upper bound on MSEE of state 1} that 
\begin{equation}
\label{eqn:KFSS trace bound}
\textrm{trace}(\tilde{\Sigma})\le\tilde{\sigma}_v^2+3+3m=(m+1)(\sigma_v^2+3). 
\end{equation}
Since adding more sensors does not increase the MSEE of the corresponding Kalman filter, we have from \eqref{eqn:KFSS trace bound} that $\textrm{trace}(\Sigma(\tilde{\mu}))\le (m+1)(\sigma_v^2+3)$, which further implies $\textrm{trace}(\Sigma(\mu^*))\le (m+1)(\sigma_v^2+3)$, where $\mu^*$ is the optimal sensor selection of the priori KFSS problem. Since $\mathcal{A}$ has approximation ratio $K$, it returns a sensor selection $\mu$ such that $\textrm{trace}(\Sigma(\mu))\le K(m+1)(\sigma_v^2+3)$.

Conversely, suppose that the answer to the $X3C$ instance is ``no''. Then, for any union of $l\le m$ ($l\in\mathbb{Z}_{\ge0}$) subsets in $\mathcal{C}$, denoted as $\mathcal{C}_l$, there exist $\kappa\ge1$ $(\kappa\in\mathbb{Z})$ elements in $D$ that are not covered by $\mathcal{C}_l$, i.e., for any $l\le m$ and $\mathcal{L}\triangleq\{i_1,\dots,i_l\}\subseteq\{1,\dots,\tau\}$, $G_{\mathcal{L}}\triangleq\left[\begin{matrix} g_{i_1}\ \cdots\ g_{i_l}\end{matrix}\right]^T$ has $\kappa\ge1$ zero columns. We then show that $\textrm{trace}(\Sigma(\mu))> K(m+1)(\sigma_v^2+3)$ for all sensor selections $\mu$ (that satisfy the budget constraint). We divide our arguments into two cases.

First, for any sensor selection $\mu_1$ that does not select the first sensor, the first column of $C(\mu_1)$ is zero (from the form of $C$ defined in Eq. \eqref{eqn: C matrix inapprox}). We then know from Lemma $\ref{Lemma:minimum trace of sigma}$(b) that $(\Sigma(\mu_1))_{11}=\frac{1}{1-\lambda_1^2}$. Hence, by our choice of $\lambda_1$, we have 
\begin{align}\nonumber
&(\Sigma(\mu_1))_{11}=\frac{Z^2}{Z-1/4}>Z\ge K(m+1)(\sigma_v^2+3)\\
\Rightarrow\ &\textrm{trace}(\Sigma(\mu_1))>K(m+1)(\sigma_v^2+3)\label{eqn:proof of thm 1 display eq 1},
\end{align}
where \eqref{eqn:proof of thm 1 display eq 1} follows from $\sum_{i = 2}^{3m+1}(\Sigma(\mu_1))_{ii} = 3m>0$ for all possible sensor selections. 

Second, consider sensor selections $\mu_2$ that select the first sensor. To proceed, we first assume that the measurement noise covariance is $V=\mathbf{0}_{(\tau+1)\times(\tau+1)}$. Denote $\textrm{supp}(\mu_2)=\{1,i_1,\dots,i_l\}$, where $l\le m$ and define $G(\mu_2)=\left[\begin{matrix}g_{i_1-1} & \cdots & g_{i_l-1}\end{matrix}\right]^T$. We then have
\begin{equation*}
C(\mu_2)=\begin{bmatrix}
1 & \varepsilon \mathbf{1}_{3m}^T\\
\mathbf{0} & G(\mu_2)
\end{bmatrix},
\label{eqn:C_mu_thm_1}
\end{equation*}
where $G(\mu_2)$ has $\kappa\ge1$ zero columns. As argued in Lemma $\ref{lemma:X3C no solution}$, there exists an orthogonal matrix $E\in\mathbb{R}^{(3m+1)\times (3m+1)}$ of the form $E = \left[\begin{smallmatrix}1 & 0 \\ 0 & N\end{smallmatrix}\right]$ such that 
\begin{equation*}
\tilde{C}(\mu_2)\triangleq C(\mu_2)E=\begin{bmatrix}
1 & \mathbf{\varepsilon\gamma} &\mathbf{\varepsilon\beta}\\
\mathbf{0} & \mathbf{0} & \tilde{G}(\mu_2)
\end{bmatrix}.
\end{equation*}
In the above expression, $\tilde{G}(\mu_2)\in\mathbb{R}^{l\times r}$ is of full column rank, where $r=\textrm{rank}(G(\mu_2))$. Furthermore, $\gamma\in\mathbb{R}^{1\times(3m-r)}$ and at least $\kappa$ of its elements are $1$'s, and $\beta\in\mathbb{R}^{1\times r}$. We then perform a similarity transformation on the system with $E$, which does not affect the trace of the steady state {\it a priori} error covariance matrix of the corresponding Kalman filter,\footnote{This can be easily verified using Eq. \eqref{eqn:DARE} as $E$ is an orthogonal matrix.} and does not change $A$, $W$ and $V$. We further perform additional elementary row operations to transform $\tilde{C}(\mu_2)$ into the matrix
\begin{equation*}
\tilde{C}'(\mu_2)=\begin{bmatrix}
1 & \mathbf{\varepsilon\gamma} &\mathbf{0}\\
\mathbf{0} & \mathbf{0} & \tilde{G}(\mu_2)
\end{bmatrix}.
\label{eqn:C matrix after row and column operations}
\end{equation*}
Since $A$ and $W$ are both diagonal, and $V=\mathbf{0}$, we can obtain from Eq. $(\ref{eqn:DARE})$ that the steady state {\it a priori} error covariance corresponding to the sensing matrix $\tilde{C}'(\mu_2)$, denoted as $\tilde{\Sigma}'(\mu_2)$, is of the form
\begin{equation*}
\tilde{\Sigma}'(\mu_2)=\begin{bmatrix}
\tilde{\Sigma}'_1(\mu_2) & \mathbf{0}\\
\mathbf{0} & \tilde{\Sigma}'_2(\mu_2)
\end{bmatrix},
\end{equation*}
where $\tilde{\Sigma}'_1(\mu_2)\in\mathbb{R}^{(3m+1-r)\times (3m+1-r)}$ satisfies
\begin{multline*}
\tilde{\Sigma}'_1(\mu_2)=A_1\tilde{\Sigma}'_1(\mu_2)A_1^T+W_1-\\A_1\tilde{\Sigma}'_1(\mu_2) \tilde{C}^T\big(\tilde{C}\tilde{\Sigma}'_1(\mu_2) \tilde{C}^T\big)^{-1}\tilde{C}\tilde{\Sigma}'_1(\mu_2) A_1^T,
\end{multline*}
where $A_1=\textrm{diag}(\lambda_1,0,\dots,0)\in\mathbb{R}^{(3m+1-r)\times(3m+1-r)}$, $\tilde{C}=[1\ \varepsilon\gamma]$ and $W_1=I_{3m+1-r}$. Denoting $\alpha^2=\varepsilon^2\lVert\gamma\rVert_2^2\ge\kappa\varepsilon^2\ge\varepsilon^2$, we then obtain from Lemma $\ref{lemma:estimation error of state 1}$(a) that
\begin{align}\nonumber
(\tilde{\Sigma}'(\mu_2))_{11}&=\frac{1+\alpha^2\lambda_1^2-\alpha^2+\sqrt{(\alpha^2-\alpha^2\lambda_1^2-1)^2+4\alpha^2}}{2}\\
&\ge\frac{1+\varepsilon^2\lambda_1^2-\varepsilon^2+\sqrt{(\varepsilon^2-\varepsilon^2\lambda_1^2-1)^2+4\varepsilon^2}}{2}.\label{eqn:priori msee first state}
\end{align}
By our choices of $\lambda_1$ and $\varepsilon$, we have the following:
\begin{equation}
\label{eqn:priori derivations}
\begin{split}
&\varepsilon^2>4Z^2(Z-1)\ \Rightarrow\ (1-\frac{Z-1/4}{Z})\varepsilon^2>Z^2-Z\\
\Rightarrow\ & \varepsilon^2>Z^2+Z\varepsilon^2\frac{Z-1/4}{Z^2}-Z\\
\Rightarrow\ & \varepsilon^2>Z^2+Z(\varepsilon^2(1-\lambda_1^2)-1)\\
\Rightarrow\ & (\varepsilon^2-\varepsilon^2\lambda_1^2-1)^2+4\varepsilon^2>\\
&(\varepsilon^2-\varepsilon^2\lambda_1^2-1)^2+4Z^2+4Z(\varepsilon^2-\varepsilon^2\lambda_1^2-1)\\
\Rightarrow\ & (\varepsilon^2-\varepsilon^2\lambda_1^2-1)^2+4\varepsilon^2>(2Z+\varepsilon^2-\varepsilon^2\lambda_1^2-1)^2\\
\Rightarrow\ & \sqrt{(\varepsilon^2-\varepsilon^2\lambda_1^2-1)^2+4\varepsilon^2}>2Z+\varepsilon^2-\varepsilon^2\lambda_1^2-1\\
\Rightarrow\ & \frac{1+\varepsilon^2\lambda_1^2-\varepsilon^2+\sqrt[]{(\varepsilon^2-\varepsilon^2\lambda_1^2-1)^2+4\varepsilon^2}}{2}>Z.
\end{split}
\end{equation}
Since $Z\ge K(m+1)(\sigma_v^2+3)$, \eqref{eqn:priori msee first state} and \eqref{eqn:priori derivations} imply $(\tilde{\Sigma}'(\mu_2))_{11}>K(m+1)(\sigma_v^2+3)$, which further implies $\textrm{trace}(\tilde{\Sigma}'(\mu_2))> K(m+1)(\sigma_v^2+3)$. Since $\textrm{trace}(\tilde{\Sigma}'(\mu_2))=\textrm{trace}(\Sigma(\mu_2))$ as argued above, we obtain that $\textrm{trace}(\Sigma(\mu_2))> K(m+1)(\sigma_v^2+3)$. We then note the fact that the MSEE of the Kalman filter with noiseless measurements is no greater than that with any noisy measurements (for fixed $A$, $W$ and $C$), when the system noise and the measurement noise are uncorrelated. Therefore, for $V=\sigma_v^2\begin{bmatrix}1 & \mathbf{0}\\ \mathbf{0} & \frac{1}{\varepsilon^2}I_{\tau}\end{bmatrix}$, we also have that $\textrm{trace}(\Sigma(\mu_2))> K(m+1)(\sigma_v^2+3)$ for all $\mu_2$. 

It then follows from the above arguments that $\textrm{trace}(\Sigma(\mu))> K(m+1)(\sigma_v^2+3)$ for all sensor selections $\mu$, which implies that algorithm $\mathcal{A}$ would also return a sensor selection $\mu$ such that $\text{trace}(\Sigma(\mu))>K(m+1)(\sigma_v^2+3)$. This completes the proof of the converse direction of the claim above. 

Hence, it is clear that algorithm $\mathcal{A}$ can be used to solve the $X3C$ problem by applying it to the above instance of the priori KFSS problem. Since $X3C$ is NP-complete, there is no polynomial-time algorithm for it if P $\neq$ NP, and we get a contradiction. This completes the proof of the theorem.\hfill\IEEEQED

\subsection*{Proof of Corollary $\ref{coro:inapprox KFSS post}$:}
We have shown in Theorem $\ref{thm:inapprox KFSS}$ that for any polynomial-time algorithm $\mathcal{A}$ for the priori KFSS problem and any $K\in\mathbb{R}_{\ge1}$, there exist instances of the priori KFSS problem such that $r_{\mathcal{A}}(\Sigma)>K$ (unless P $=$ NP). Suppose that there exists a polynomial-time constant-factor  approximation algorithm $\mathcal{A}'$ for the posteriori KFSS problem, i.e., $\exists K'\in\mathbb{R}_{\ge1}$ such that $r_{\mathcal{A}'}(\Sigma^*)\le K'$ for all instances of the posteriori KFSS problem, where $r_{\mathcal{A}'}(\Sigma^*)$ is defined in Eq. \eqref{eqn:performance ratio of algorithm post}. We consider an instance of the priori KFSS problem constructed in the proof of Theorem $\ref{thm:inapprox KFSS}$. We then set the instance of the posteriori KFSS problem to be the same as the constructed instance of the priori KFSS problem. Since $A=\textrm{diag}(\lambda_1,0,\dots,0)\in\mathbb{R}^{(3m+1)\times(3m+1)}$ and $W=I_{3m+1}$, where $0<\lambda_1<1$, we have from Eq. \eqref{eqn:couple priori and posteriori 1} that
\begin{equation}
\label{eqn:coupled priori and posteriori state 1}
(\Sigma(\mu))_{11}=\lambda_1^2(\Sigma^*(\mu))_{11}+1,\forall\mu.
\end{equation}
Since we know from Lemma $\ref{Lemma:minimum trace of sigma}$(a) that $(\Sigma(\mu))_{ii}=1$, $\forall i\in\{2,\dots,3m+1\}$ and $\forall\mu$, it then follows from Eq. \eqref{eqn:coupled priori and posteriori state 1} that
\begin{equation}
\label{eqn:coupled trace}
\textrm{trace}(\Sigma(\mu))=\lambda_1^2(\Sigma^*(\mu))_{11}+3m+1, \forall\mu.
\end{equation}
We also know from Lemma $\ref{Lemma:minimum trace of sigma}$(a) that $0\le(\Sigma^*(\mu))_{ii}\le1$, $\forall i\in\{2,\dots,3m+1\}$ and $\forall\mu$, which implies that
\begin{equation}
\label{eqn:upperbound on posteriori trace 1}
\textrm{trace}(\Sigma^*(\mu))\le(\Sigma^*(\mu))_{11}+3m,\forall\mu.
\end{equation}
We then obtain from Eqs. \eqref{eqn:coupled trace}-\eqref{eqn:upperbound on posteriori trace 1} that 
\begin{align}\nonumber
\textrm{trace}(\Sigma^*(\mu))&\le \frac{3m\lambda_1^2+\textrm{trace}(\Sigma(\mu))-3m-1}{\lambda_1^2}\\
&\le \frac{\textrm{trace}(\Sigma(\mu))}{\lambda_1^2}, \forall\mu\label{eqn:upperbound on posteriori trace 2},
\end{align}
where the second inequality follows from the fact that $0<\lambda_1<1$.
Denote the optimal sensor selections of the priori and the posteriori KFSS problems as $\mu_1^*$ and $\mu^*_2$, respectively. Denote the sensor selection returned by algorithm $\mathcal{A}'$ as $\mu'$. Note that $\Sigma_{opt}=\Sigma(\mu^*_1)$ and $\Sigma^*_{opt}=\Sigma^*(\mu_2^*)$ and $\Sigma^*_{\mathcal{A}'}=\Sigma^*(\mu')$. We then have the following:
\begin{align}\nonumber
&\text{trace}(\Sigma^*_{\mathcal{A}'})\le K'\text{trace}(\Sigma^*_{opt})\\\nonumber
\Rightarrow\ &(\Sigma^*(\mu'))_{11}+\sum_{i=2}^{3m+1} (\Sigma^*(\mu'))_{ii}\le K' \text{trace}(\Sigma^*(\mu^*_2))\\
\Rightarrow\ &\frac{(\Sigma(\mu'))_{11}-1}{\lambda_1^2}\le  K'\text{trace}(\Sigma^*(\mu^*_2))\le K'\text{trace}(\Sigma^*(\mu^*_1))\label{eqn:coro KFSS post 1}\\
\Rightarrow\ &(\Sigma(\mu'))_{11}-1\le K'\textrm{trace}(\Sigma(\mu^*_1))\label{eqn:coro KFSS post 1_1} \\
\Rightarrow\ &\textrm{trace}(\Sigma(\mu'))\le K'\textrm{trace}(\Sigma(\mu^*_1))+3m+1\label{eqn:coro KFSS post 2}\\
\Rightarrow\ &\frac{\textrm{trace}(\Sigma(\mu'))}{\textrm{trace}(\Sigma(\mu_1^*))}\le K'+\frac{3m+1}{\textrm{trace}(\Sigma(\mu^*_1))}\le K'+1 \label{eqn:coro KFSS post 3},
\end{align}
where the first inequality in \eqref{eqn:coro KFSS post 1} follows from Eq. \eqref{eqn:coupled priori and posteriori state 1} and $(\Sigma^*(\mu'))_{ii}\ge0,\forall i$ (from Lemma $\ref{Lemma:minimum trace of sigma}$(a)), the second inequality in \eqref{eqn:coro KFSS post 1} follows from the fact that $\mu_2^*$ is the optimal sensor selection for the posteriori KFSS problem,  \eqref{eqn:coro KFSS post 1_1} follows from \eqref{eqn:upperbound on posteriori trace 2}, \eqref{eqn:coro KFSS post 2} follows from the fact that $\sum_{i=2}^{3m+1}(\Sigma(\mu'))_{ii}=3m$ (from Lemma $\ref{Lemma:minimum trace of sigma}$(a)), and the second inequality in \eqref{eqn:coro KFSS post 3} uses the fact that $\textrm{trace}(\Sigma(\mu^*_1))\ge 3m+1$ (from Lemma $\ref{Lemma:minimum trace of sigma}$(a)). Thus, we have from \eqref{eqn:coro KFSS post 3} that $r_{\mathcal{A'}}(\Sigma)\le K'+1$, which contradicts the fact that the priori KFSS problem cannot have a polynomial-time constant-factor  approximation algorithm for instances of the given form, and completes the proof of the corollary.
\hfill\IEEEQED

\section*{Appendix C}
\subsection*{Proof of Theorem $\ref{thm:inapprox KFSA}$:}
Assume that there exists such a polynomial-time constant-factor approximation algorithm $\mathcal{A}$, i.e., $\exists K\in\mathbb{R}_{\ge1}$ such that $r_{\mathcal{A}}(\tilde{\Sigma})\le K$ for all instances of the priori KFSA problem, where $r_{\mathcal{A}}(\tilde{\Sigma})$ is defined in Eq. $\eqref{eqn:performance ratio of KFSA algorithms}$. We will show that $\mathcal{A}$ can be used to solve the $X3C$ problem, leading to a contradiction.

Consider any instance of the $X3C$ problem to be a finite set $D=\{d_1,\cdots,d_{3m}\}$ and a collection $\mathcal{C}=\{c_1,\dots,c_{\tau}\}$ of $3$-element subsets of $D$, where $\tau\ge m$. Recall in the proof of Theorem $\ref{thm:inapprox KFSS}$ that we use a column vector $g_i \in \mathbb{R}^{3m}$ to encode which elements of $D$ are contained in $c_i$, where $(g_i)_j = 1$ if $d_j\in D$ is in $c_i$, and $(g_i)_j = 0$ otherwise, for $i \in \{1, 2, \ldots, \tau\}$ and $j \in \{1, 2, \ldots, 3m\}$. The matrix $G\in\mathbb{R}^{\tau\times 3m}$ was defined in the proof of Theorem $\ref{thm:inapprox KFSS}$ as $G=\left[\begin{matrix}g_1 & \cdots & g_{\tau}\end{matrix}\right]^T$. In this proof, we will make use of the matrix $F\triangleq G^T$; note that each column of $F$ contains exactly three $1$'s.

Given the above instance of the $X3C$ problem, we then construct an instance of the priori KFSA as follows. Denote $Z=\lceil K\rceil(\tau+2)(\delta_v^2+1)$, where we set $\delta_v=1$. Define the system dynamics matrix as $A=\text{diag}(\lambda_1,0,\dots,0)\in\mathbb{R}^{(\tau+1)\times(\tau+1)}$, where $\lambda_1=\frac{Z-1/2}{Z}$. Note that $Z\in\mathbb{Z}_{>1}$ and $0<\lambda_1<1$. The set $\mathcal{Q}$ consists of $3m+\tau$ sensors with collective measurement matrix
\begin{equation}
C=\begin{bmatrix}
\mathbf{1}_{3m} & \rho F\\
\mathbf{0} & I_{\tau}
\end{bmatrix},
\label{eqn:C matrix in Theorem 2}
\end{equation}
where $F$ is defined above and $I_{\tau}$ is used to encode the collection $\mathcal{C}$, i.e., $\mathbf{e}_{j}$ represents $c_j\in \mathcal{C}$ for all $j\in\{1,2,\dots,\tau\}$. The constant $\rho$ is chosen as $\rho=2Z\Bigl\lceil\sqrt{m(Z-1)}\Bigr\rceil+1$. The system noise covariance matrix is set to be $W=I_{\tau+1}$. The measurement noise covariance is set as $V=\delta_v^2\begin{bmatrix}I_{3m} & \mathbf{0}\\ \mathbf{0} & \frac{1}{\rho^2}I_{\tau}\end{bmatrix}$. The sensor attack cost vector is set as $\omega=\mathbf{1}_{3m+\tau}$, and the sensor attack budget is set as $\Omega=m$. Note that the sensor attack vector is given by $\nu\in\{0,1\}^{3m+\tau}$.  For the above construction, since the only nonzero eigenvalue of $A$ is $\lambda_1$, we know from Lemma $\ref{Lemma:minimum trace of sigma}$(a) that $\sum_{i=2}^{\tau+1}(\Sigma(\nu^c))_{ii}=\sum_{i=2}^{\tau+1}W_{ii}=\tau$ for all $\nu$. 

We claim that algorithm $\mathcal{A}$ will return a sensor attack vector $\nu$ such that $\textrm{trace}(\Sigma(\nu^c))>(\tau+2)(\delta_v^2+1)$ if and only if the answer to the $X3C$ problem is ``yes''. 

We prove the above claim as follows. Suppose that the answer to the $X3C$ problem is ``yes''. Similarly to the proof of Theorem $\ref{thm:inapprox KFSS}$, we first assume that $V=\mathbf{0}_{(3m+\tau)\times(3m+\tau)}$. Denote an exact cover as $\mathcal{C}'=\{c_{j_1},\dots, c_{j_m}\}$, where $\{j_{1},\dots, j_{m}\}\subseteq\{1,2,\dots,\tau\}$. Define $\tilde{\nu}$ to be the sensor attack such that $\textrm{supp}(\tilde{\nu})=\{3m+j_1,\dots,3m+j_m\}$. We then renumber the states of the system from state $2$ to state $\tau$ such that for all $i\in\{1,2,\dots, m\}$, the columns of the submatrix $I_{\tau}$ of $C$ in Eq. \eqref{eqn:C matrix in Theorem 2} representing $c_{j_{i}}$ in $\mathcal{C}'$, i.e., the columns of $I_{\tau}$ that correspond to $\text{supp}(\tilde{\nu})$, come first. Note that renumbering the states does not change the trace of the steady state {\it a priori} error covariance of the corresponding Kalman filter. We then have from Eq. $\eqref{eqn:C matrix in Theorem 2}$ that
\begin{equation}
C(\tilde{\nu}^c)=\begin{bmatrix}
\mathbf{1}_{3m} & \rho F_1 & \rho F_2\\
\mathbf{0} & \mathbf{0} & I_{\tau-m}
\end{bmatrix},
\label{eqn:tildeC_nu in Theorem 2}
\end{equation}
where $F_1\in\mathbb{R}^{3m\times m}$ and $F_2\in\mathbb{R}^{3m\times(\tau-m)}$ satisfy $F=\begin{bmatrix} F_1 & F_2\end{bmatrix}$, and $I_{\tau-m}$ is the submatrix of $I_{\tau}$ that corresponds to $\textrm{supp}(\tilde{\nu}^c)\cap\{3m+1,\dots,3m+\tau\}$, i.e., the elements of $\mathcal{C}$ that are not in $\mathcal{C}'$.\footnote{Note that if the submatrix of $I_{\tau}$ corresponding to $\text{supp}(\tilde{\nu}^c)\cap\{3m+1,\dots,3m+\tau\}$ is not identity, we can always permute the rows of $C(\tilde{\nu}^c)$ to make it identity.} Since the sensor attack $\tilde{\nu}$ targets the rows of $C$ that correspond to the elements of the exact cover $\mathcal{C}'$ for $D$, we have that $F_1$, after some row permutations of $C(\tilde{\nu}^c)$, is given by $F_1=\left[\begin{matrix}\mathbf{e}_1^T & \mathbf{e}_1^T & \mathbf{e}_1^T\cdots & \mathbf{e}_m^T & \mathbf{e}_m^T & \mathbf{e}_m^T \end{matrix}\right]^T$. We perform additional elementary row operations and merge identical rows (which does not change the steady state {\it a priori} error covariance matrix of the corresponding Kalman filter) to transform $C(\tilde{\nu}^c)$ into the matrix
\begin{equation}
\tilde{C}(\tilde{\nu}^c)=\begin{bmatrix}
\mathbf{1}_m & \rho I_m & \mathbf{0}\\
\mathbf{0} & \mathbf{0} & I_{\tau-m}
\end{bmatrix}.
\label{eqn:tildeC_nu_2 in Theorem 2}
\end{equation}
Since $A$ and $W$ are both diagonal, and $V=\mathbf{0}$, we can obtain from Eq. $(\ref{eqn:DARE})$ that the steady state {\it a priori} error covariance corresponding to $\tilde{C}(\tilde{\nu}^c)$, denoted as $\tilde{\Sigma}(\tilde{\nu}^c)$, is of the form
\begin{equation*}
\tilde{\Sigma}(\tilde{\nu}^c)=\begin{bmatrix}
\tilde{\Sigma}_1(\tilde{\nu}^c) & \mathbf{0}\\
\mathbf{0} & \tilde{\Sigma}_2(\tilde{\nu}^c)
\end{bmatrix},
\end{equation*}
where $\tilde{\Sigma}_1(\tilde{\nu}^c)\in\mathbb{R}^{(m+1)\times (m+1)}$ satisfies
\begin{multline*}
\tilde{\Sigma}_1(\tilde{\nu}^c)=A_1\tilde{\Sigma}_1(\tilde{\nu}^c) A_1^T+W_1-\\A_1\tilde{\Sigma}_1(\tilde{\nu}^c) \tilde{C}^T\big(\tilde{C}\tilde{\Sigma}_1(\tilde{\nu}^c) \tilde{C}^T\big)^{-1}\tilde{C}\tilde{\Sigma}_1(\tilde{\nu}^c) A_1^T,
\end{multline*}
where $A_1=\textrm{diag}(\lambda_1,0,\dots,0)\in\mathbb{R}^{(m+1)\times(m+1)}$, $\tilde{C}=\left[\begin{matrix}\mathbf{1}_m & \rho I_m \end{matrix}\right]$ and $W_1=I_{m+1}$.  We then know from Lemma $\ref{lemma:estimation error of state 1}$(b) that $(\Sigma(\tilde{\nu}^c))_{11}=(\tilde{\Sigma}(\tilde{\nu}^c))_{11}$ satisfies
\begin{equation}
(\Sigma(\tilde{\nu}^c))_{11}=\frac{\lambda_1^2\rho^2+m-\rho^2+\sqrt{(\rho^2-\lambda_1^2\rho^2-m)^2+4m\rho^2}}{2m}.
\label{eqn:sigma11 after attack}
\end{equation}
By our choices of $\lambda_1$ and $\rho$, we have 
\begin{equation}
\label{eqn:priori derivations KFSA}
\begin{split}
&\rho^2>4Z^2m(Z-1)\Rightarrow  (1-\frac{Z-1/4}{Z})\rho^2>Z^2m-Zm\\
\Rightarrow\ & \rho^2>mZ^2+Z\rho^2\frac{Z-1/4}{Z^2}-Zm\\
\Rightarrow\ & 4m\rho^2>4m^2Z^2+4mZ(\rho^2(1-\lambda_1^2)-m)\\
\Rightarrow\ & (\rho^2-\lambda_1^2\rho^2-m)^2+4m\rho^2\\
&>4m^2Z^2+4mZ(\rho^2-\lambda_1^2\rho^2-m)+(\rho^2-\lambda_1^2\rho^2-m)^2\\
\Rightarrow\ & (\rho^2-\lambda_1^2\rho^2-m)^2+4m\rho^2>(2mZ+\rho^2-\lambda_1^2\rho^2-m)^2\\
\Rightarrow\ & \sqrt{ (\rho^2-\lambda_1^2\rho^2-m)^2+4m\rho^2}>2mZ+\rho^2-\lambda_1^2\rho^2-m\\
\Rightarrow\ & \frac{\lambda_1^2\rho^2+m-\rho^2+\sqrt{(\rho^2-\lambda_1^2\rho^2-m)^2+4m\rho^2}}{2m}>Z.
\end{split}
\end{equation}
Noting that $Z\ge K(\tau+2)(\delta_v^2+1)$, we then know from \eqref{eqn:sigma11 after attack} and \eqref{eqn:priori derivations KFSA} that $(\Sigma(\tilde{\nu}^c))_{11}>K(\tau+2)(\delta_v^2+1)$, which further implies that $\text{trace}(\Sigma(\tilde{\nu}^c))>K(\tau+2)(\delta_v^2+1)$.  Following the same arguments as those in the proof of Theorem $\ref{thm:inapprox KFSS}$, we have that for $V=\delta_v^2\begin{bmatrix}I_{3m} & \mathbf{0}\\ \mathbf{0} & \frac{1}{\rho^2}I_{\tau}\end{bmatrix}$, $\text{trace}(\Sigma(\tilde{\nu}^c))>K(\tau+2)(\delta_v^2+1)$ also holds, which implies $\text{trace}(\Sigma(\nu^{*c}))>K(\tau+2)(\delta_v^2+1)$, where $\nu^*$ is the optimal sensor attack for the priori KFSA problem. Since algorithm $\mathcal{A}$ has approximation ratio $K$, it would return a sensor attack $\nu$ such that $\text{trace}(\Sigma(\nu^c))>(\tau+2)(\delta_v^2+1)$. 

Conversely, suppose the answer to the $X3C$ problem is ``no''. For any union of $l\le m$ ($l\in\mathbb{Z}_{\ge0}$) subsets in $\mathcal{C}$, denoted as $\mathcal{C}_l$, there exists at least one element in $D$ that is not covered by $\mathcal{C}_{l}$. We then show that $\text{trace}(\Sigma(\nu^c))\le(\tau+2)(\delta_v^2+1)$ for all sensor attacks $\nu$ (that satisfy the budget constraint). We split our discussion into three cases.

First, consider any sensor attack $\nu_1$ that targets $l$ sensors merely from $C_1$ to $C_{3m}$ in Eq. $\eqref{eqn:C matrix in Theorem 2}$, i.e., $|\textrm{supp}(\nu_1)|=l$ and $\textrm{supp}(\nu_1)\subseteq\{1,\dots,3m\}$, where $l\le m$. We then obtain 
\begin{equation*}
C(\nu_1^c)=\begin{bmatrix}
\mathbf{1}_{3m-l} & \rho F(\nu_1^c)\\
\mathbf{0} & I_{\tau}
\end{bmatrix},
\end{equation*}
where $F(\nu^c_1)\in\mathbb{R}^{(3m-l)\times\tau}$ is defined to be the submatrix of $F$ that corresponds to $\textrm{supp}(\nu_1^c)\cap\{1,\dots,3m\}$, i.e., the rows of $F$ that are left over by $\nu_1$. We perform elementary row operations to transform $C(\nu_1^c)$ into
\begin{equation}
\tilde{C}(\nu_1^c)\triangleq \Psi C(\nu_1^c)=\begin{bmatrix}
\mathbf{1}_{3m-l} & \mathbf{0}\\
\mathbf{0} & I_{\tau}
\end{bmatrix}
\label{eqn:C after row operation}
\end{equation}
with the corresponding measurement noise covariance 
\begin{equation}
\label{eqn:covariance after row ops}
\begin{split}
\tilde{V}(\nu_1^c)&\triangleq \Psi V(\nu_1^c)\Psi^T\\
&=\delta_v^2\begin{bmatrix}I_{3m-l}+F(\nu_1^c)(F(\nu_1^c))^T & -\frac{1}{\rho}F(\nu_1^c)\\ -\frac{1}{\rho}(F(\nu_1^c))^T & \frac{1}{\rho^2}I_{\tau} \end{bmatrix},
\end{split}
\end{equation} 
where $\Psi=\begin{bmatrix}I_{3m-l} & -\rho F(\nu_1^c)\\ \mathbf{0} & I_{\tau}\end{bmatrix}$. Since there are at most $\tau$ nonzero elements (which are all $1$'s) in the first row of $F(\nu_1^c)$, it follows that $(F(\nu_1^c)(F(\nu_1^c))^T)_{11}$ (i.e., the element in the first row and first column of the matrix $F(\nu_1^c)(F(\nu_1^c))^T$) is at most $\tau$. We then have from Eq. \eqref{eqn:covariance after row ops} that $(\tilde{V}(\nu_1^c))_{11}$, denoted as $\tilde{\delta}_v^2(\nu_1^c)$, satisfies
\begin{equation}
\label{eqn:bound on delta_1}
\tilde{\delta}_v^2(\nu_1^c)\le(\tau+1)\delta_v^2. 
\end{equation}

Second, consider any sensor attack $\nu_2$ that targets $l$ sensors merely from $C_{3m+1}$ to $C_{3m+\tau}$ in Eq. $\eqref{eqn:C matrix in Theorem 2}$, i.e.,  $|\textrm{supp}(\nu_2)|=l$ and $\text{supp}(\nu_2)\subseteq\{3m+1,\dots,3m+\tau\}$, where $l\le m$. Via similar arguments to those for obtaining Eqs. \eqref{eqn:tildeC_nu in Theorem 2}, \eqref{eqn:C after row operation} and \eqref{eqn:covariance after row ops}, we can perform elementary row operations to transform 
\begin{equation*}
C(\nu_2^c)=\begin{bmatrix}
\mathbf{1}_{3m} & \rho F'_1 & \rho F'_2\\
\mathbf{0} & \mathbf{0} & I_{\tau-l}
\end{bmatrix}
\end{equation*}
 into
\begin{equation*}
\tilde{C}(\nu_2^c)=\begin{bmatrix}
\mathbf{1}_{3m} & \rho F'_1 & \mathbf{0}\\
\mathbf{0} & \mathbf{0} & I_{\tau-l}
\end{bmatrix}
\end{equation*}
with the corresponding measurement noise covariance $\tilde{V}(\nu_2^c)=\begin{bmatrix}\tilde{\delta}_v^2(\nu_2^c) & *\\ * & *\end{bmatrix}$, where 
\begin{equation}
\label{eqn:bound on delta_2}
\tilde{\delta}_v^2(\nu_2^c)\le(\tau-l+1)\delta_v^2.
\end{equation}
Note that $F'_1\in\mathbb{R}^{3m\times l}$ and $F'_2\in\mathbb{R}^{3m\times(\tau-l)}$ satisfy $F=\begin{bmatrix} F'_1 & F'_2\end{bmatrix}$. Recall that for any union of $l\le m$ subsets in $\mathcal{C}$, denoted as $\mathcal{C}_l$, there exists at least one element in $D$ that is not covered by $\mathcal{C}_l$. We can then assume without loss of generality that one such element is $d_1$, which implies that the first row of $F'_1$ is zero. 

Third, consider any sensor attack $\nu_3$ that targets sensors from  both $C_1$ to $C_{3m}$ and $C_{3m+1}$ to $C_{3m+\tau}$ in Eq. $\eqref{eqn:C matrix in Theorem 2}$. Suppose that the attack $\nu_3$ attacks $l_1$ sensors from  $C_1$ to $C_{3m}$ and $l_2$ sensors from $C_{3m+1}$ to $C_{3m+\tau}$, i.e., $\text{supp}(\nu_3)=\{j'_1,\dots,j'_{l_1},3m+j''_1,\dots,3m+j''_{l_2}\}\subseteq\{1,2,\dots,3m+\tau\}$, where $l_1,l_2\in\mathbb{Z}_{\ge1}$, $l_1+l_2=l\le m$, $\{j'_1,\dots,j'_{l_1}\}\subseteq\{1,\dots,3m\}$ and $\{j''_1,\dots,j''_{l_2}\}\subseteq\{1,\dots,\tau\}$. By similar arguments to those above, we can perform elementary row operations to transform 
\begin{equation*}
C(\nu_3^c)=\begin{bmatrix}
\mathbf{1}_{3m-l_1} & \rho F_1(\nu_3^c) & \rho F_2(\nu_3^c)\\
\mathbf{0} & \mathbf{0} & I_{\tau-l_2}
\end{bmatrix}
\end{equation*}
into
\begin{equation*}
\tilde{C}(\nu_3^c)=\begin{bmatrix}
\mathbf{1}_{3m-l_1} & \rho F_1(\nu_3^c) & 0\\
\mathbf{0} & \mathbf{0} & I_{\tau-l_2}
\end{bmatrix},
\end{equation*}
where $F_1(\nu_3^c)\in\mathbb{R}^{(3m-l_1)\times l_2}$ and $F_2(\nu_3^c)\in\mathbb{R}^{(3m-l_1)\times(\tau- l_2)}$ satisfy $F(\nu_3^c)=\begin{bmatrix}F_1(\nu_3^c) & F_2(\nu_3^c)\end{bmatrix}$ with $F(\nu_3^c)$ defined in the same way as $F(\nu_1^c)$. Moreover, the measurement noise covariance corresponding to $\tilde{C}(\nu_3^c)$ is given by $\tilde{V}(\nu_3^c)=\begin{bmatrix}\tilde{\delta}_v^2(\nu_3^c) & *\\ * & *\end{bmatrix}$, where
\begin{equation}
\label{eqn:bound on delta_3}
 \tilde{\delta}_v^2(\nu_3^c)\le(\tau-l_2+1)\delta_v^2.
 \end{equation}
Since any $l_2$ subsets in $\mathcal{C}$ can cover at most $3l_2$ elements in $D$, there are at least $3m-3l_2$ elements in $D$ that are not covered by the $l_2$ subsets in $\mathcal{C}$. Also note that 
\begin{align*}
3m-3l_2-l_1=3m-2l_2-l=2(m-l_2)+m-l>0,
\end{align*}
where the last inequality follows from the facts that $l_1+l_2=l\le m$ and $l_1,l_2\in\mathbb{Z}_{\ge1}$. Hence, by attacking $l_1$ sensors from  $C_1$ to $C_{3m}$ and $l_2$ sensors from $C_{3m+1}$ to $C_{3m+\tau}$, we have at least $3m-3l_2-l_1>0$ row(s) of $F_1(\nu_3^c)$ that are zero. Again, we can assume without loss of generality that the first row of $F_1(\nu_3^c)$ is zero. 

In summary, for any sensor attack $\nu$, we let $\hat{\Sigma}(\nu_i^c)$ denote the steady state {\it a priori} error covariance obtained from measurement matrix $(\tilde{C}(\nu_i^c))_1=\mathbf{e}_1$ with measurement noise variance $\tilde{\delta}_v^2(\nu_i^c)$ (which corresponds to the first sensor in $\tilde{C}(\nu_i^c)$), $\forall i\in\{1,2,3\}$, where $\nu_1$, $\nu_2$ and $\nu_3$ are given as above. Following similar arguments to those for \eqref{eqn:thm1 upper bound on MSEE of state 1}, we have $(\hat{\Sigma}(\nu_i^c))_{11}\le\tilde{\delta}_v^2(\nu_i^c)+2$, $\forall i\in\{1,2,3\}$. Since $\sum_{i=2}^{\tau+1}(\hat{\Sigma}(\nu^c_i))_{ii}=\sum_{i=2}^{\tau+1}W_{ii}=\tau$ holds for all $i\in\{1,2,3\}$ via similar arguments to those above, we obtain that 
\begin{equation}
\label{eqn:KFSA trace bound}
\textrm{trace}(\hat{\Sigma}(\nu_i^c))\le\tilde{\delta}_v^2(\nu_i^c)+2+\tau, \forall i\in\{1,2,3\}. 
\end{equation}
Again note that adding more sensors does not increase the MSEE of the corresponding Kalman filter, and the above operations performed on the sensing matrix $C$ do not change the trace of the steady state {\it a priori} error covariance of the corresponding Kalman filter as well. We then have from Eqs. \eqref{eqn:bound on delta_1}-\eqref{eqn:KFSA trace bound} that $\textrm{trace}(\Sigma(\nu^c))\le(\tau+1)\delta_v^2+2+\tau\le(\tau+2)(\delta_v^2+1)$ for all $\nu$. It follows that algorithm $\mathcal{A}$ would also return a sensor attack $\nu$ such that $\text{trace}(\Sigma(\nu^c))\le(\tau+2)(\delta_v^2+1)$. This proves the converse direction of the claim above. 

Therefore, we know that $\mathcal{A}$ can be used to solve the $X3C$ problem by applying it to the above instance of the priori KFSA problem. Since $X3C$ is NP-complete, there is no polynomial-time algorithm for it if P $\neq$ NP, yielding a contradiction. This completes the proof of the theorem.\hfill\IEEEQED

\subsection*{Proof of Corollary $\ref{coro:inapprox KFSA post}$:}
Note that the $A$ and $W$ matrices for the instance of KFSA that we constructed in the proof of Theorem $\ref{thm:inapprox KFSA}$ are the same as those for the instance of KFSS that we constructed in the proof of Theorem $\ref{thm:inapprox KFSS}$. We then follow the same arguments as those in the proof of Corollary $\ref{coro:inapprox KFSS post}$. Suppose that there exists a polynomial-time constant-factor  approximation algorithm $\mathcal{A}'$ for the posteriori KFSA problem, i.e., $\exists K'\in\mathbb{R}_{\ge1}$ such that $r_{\mathcal{A}'}(\tilde{\Sigma}^*)\le K'$ for all instances of the posteriori KFSA problem, where $r_{\mathcal{A}'}(\tilde{\Sigma}^*)$ is defined in Eq. \eqref{eqn:performance ratio of post KFSA algorithms}. We consider an instance of the priori KFSA problem as constructed in the proof of Theorem $\ref{thm:inapprox KFSA}$. We then set the instance of the posteriori KFSA problem to be the same as the constructed instance of the priori KFSA problem. Denote the optimal sensor attacks of the priori and the posteriori KFSA problems as $\nu_1^*$ and $\nu^*_2$, respectively. Denote the sensor attack returned by algorithm $\mathcal{A}'$ as $\nu'$. Note that $\tilde{\Sigma}_{opt}=\Sigma(\nu_1^{*c})$, $\tilde{\Sigma}^*_{opt}=\Sigma^*(\nu_2^{*c})$ and $\tilde{\Sigma}^*_{\mathcal{A}'}=\Sigma^*(\nu'^c)$. Also note that $\text{trace}(\Sigma^*(\nu_1^{*c}))\le \text{trace}(\Sigma^*(\nu^{*c}_2))$, since $\nu^*_2$ is the optimal sensor attack for the posteriori KFSA problem. We then have the following:
\begin{align*}\nonumber
&\text{trace}(\Sigma^*(\nu_1^{*c}))\le \text{trace}(\Sigma^*(\nu^{*c}_2))\le K'\text{trace}(\tilde{\Sigma}^*_{\mathcal{A}'})\\
\Rightarrow\ &(\Sigma^*(\nu_1^{*c}))_{11}+\sum_{i=2}^{3m+1} (\Sigma^*(\nu^{*c}_1))_{ii}\le K' \text{trace}(\Sigma^*(\nu'^c))\\
\Rightarrow\  &\frac{(\Sigma(\nu_1^{*c}))_{11}-1}{\lambda_1^2}\le K'\text{trace}(\Sigma^*(\nu'^c))\\
\Rightarrow\  &(\Sigma(\nu_1^{*c}))_{11}-1\le K'\textrm{trace}(\Sigma(\nu'^c)) \\
\Rightarrow\ &\textrm{trace}(\Sigma(\nu_1^{*c}))\le K'\textrm{trace}(\Sigma(\nu'^c))+3m+1\\
\Rightarrow\  &\frac{\textrm{trace}(\Sigma(\nu_1^{*c}))}{\textrm{trace}(\Sigma(\nu'^c))}\le K'+\frac{3m+1}{\textrm{trace}(\Sigma(\nu'^c))}\le K'+1,
\end{align*}
which implies $r_{\mathcal{A}'}(\tilde{\Sigma})\le K'+1$, and yields a contradiction with the fact that the priori KFSA problem cannot have a polynomial-time constant-factor  approximation algorithm for the instances of the form given as above. This completes the proof of the corollary.
\hfill\IEEEQED

\section*{Appendix D}
\subsection*{Proof of Theorem $\ref{thm:ratio of special system KFSS}$:}
We first prove that Algorithm \ref{algorithm:greedy KFSS} for the priori KFSS problem selects sensor $2$ and sensor $3$ in its first and second iterations, respectively. Since the only nonzero eigenvalue of $A$ is $\lambda_1$, we know from Lemma $\ref{Lemma:minimum trace of sigma}$(a) that $(\Sigma(\mu))_{22}=1$ and $(\Sigma(\mu))_{33}=1$, $\forall\mu$, which implies that $(\Sigma_{gre})_{22}=1$ and $(\Sigma_{gre})_{33}=1$. Hence, we focus on determining $(\Sigma_{gre})_{11}$.

Denoting $\mu_1=[1\ 0\ 0]^T$ and $\mu_2=[0\ 1\ 0]^T$, we have $C(\mu_1)=[1\ h\ h]$ and $C(\mu_2)=[1\ 0\ h]$. Using the result in Lemma $\ref{lemma:estimation error of state 1}$(a), we obtain that $\sigma_1\triangleq(\Sigma(\mu_1))_{11}$ and $\sigma_2\triangleq(\Sigma(\mu_2))_{11}$ satisfy
\begin{equation*}
\sigma_1=\frac{2}{\sqrt[]{(1-\lambda_1^2-\frac{1}{2h^2})^2+\frac{2}{h^2}}+1-\lambda_1^2-\frac{1}{2h^2}},
\end{equation*}
and
\begin{equation*}
\sigma_2=\frac{2}{\sqrt[]{(1-\lambda_1^2-\frac{1}{h^2})^2+\frac{4}{h^2}}+1-\lambda_1^2-\frac{1}{h^2}},
\end{equation*}
respectively. Similarly, denoting $\mu_3=[0\ 0\ 1]^T$, we obtain $C(\mu_3)=[0\ 1\ 1]$. Since the first column of $C(\mu_3)$ is zero, we know from Lemma $\ref{Lemma:minimum trace of sigma}$(b) that $\sigma_3 \triangleq (\Sigma(\mu_3))_{11}=\frac{1}{1-\lambda_1^2}$.  If we view $\sigma_2$ as a function of $h^2$, denoted as $\sigma(h^2)$, we have $\sigma_1=\sigma(2h^2)$. Since we know from Lemma $\ref{lemma:estimation error of state 1}$(a) that $\sigma(h^2)$ is a strictly increasing function of $h^2\in\mathbb{R}_{>0}$ and upper bounded by $\frac{1}{1-\lambda_1^2}$, we obtain $\sigma_2<\sigma_1<\sigma_3$, which implies that the greedy algorithm selects sensor $2$ in its first iteration.

Denote $\mu_{12}=[1\ 1\ 0]^T$. We have $C(\mu_{12})=
\left[\begin{smallmatrix}
    1 & h & h \\
    1 & 0 & h \\
\end{smallmatrix}\right]$, on which we perform elementary row operations and obtain $\tilde{C}(\mu_{12})=
\left[\begin{smallmatrix}
    0 & h & 0 \\
    1 & 0 & h \\
\end{smallmatrix}\right]$. By direct computation from Eq. $\eqref{eqn:DARE}$, we obtain $(\Sigma(\mu_{12}))_{11}=\sigma_2$. Moreover, we denote $\mu_{23}=[0\ 1\ 1]^T$ and obtain $C(\mu_{23})=
\left[\begin{smallmatrix}
    1 & 0 & h \\
    0 & 1 & 1 \\
\end{smallmatrix}\right]$. By direct computation from Eq. $\eqref{eqn:DARE}$, we have $(\Sigma(\mu_{23}))_{11}$, denoted as $\sigma_{23}$, to be 
\begin{equation*}
\sigma_{23}=\frac{2}{\sqrt[]{(1-\lambda_1^2-\frac{2}{h^2})^2+\frac{8}{h^2}}+1-\lambda_1^2-\frac{2}{h^2}}.
\label{eqn:solution for sigma23}
\end{equation*}
Similarly to the argument above, we have $\sigma_{12}=\sigma(h^2)$ and $\sigma_{23}=\sigma(\frac{h^2}{2})$, where $\sigma(\frac{h^2}{2})<\sigma(h^2)$, which implies that the greedy algorithm selects sensor $3$ in its second iteration. Hence, we have $\textrm{trace}(\Sigma_{gre})=\sigma_{23}+2$.

Furthermore, it is easy to see that the optimal sensor selection (for the priori KFSS instance) is $\mu=[1\ 0\ 1]^T$, denoted as $\mu_{13}$. Since if $\mu=\mu_{13}$, then $\mathbf{e}_1\in\textrm{rowspace}(C(\mu))$ and thus we know from Lemma $\ref{Lemma:minimum trace of sigma}$(a) and (c) that $\textrm{trace}(\Sigma(\mu))=3=\textrm{trace}(W)$, which is also the minimum value of $\textrm{trace}(\Sigma(\mu))$ among all possible sensor selections $\mu$. Combining the results above and taking the limit as $h\to\infty$ lead to Eq. $(\ref{eqn:limit ratio KFSS})$.

We next prove that the greedy algorithm defined in Algorithm \ref{algorithm:greedy KFSS} for the posteriori KFSS problem selects sensor $2$ and sensor $3$ in its first and second iterations, respectively. Note that it is easy to obtain from Eq. \eqref{eqn:couple priori and posteriori 1} that $\Sigma(\mu)$ is of the form $\Sigma(\mu)=\text{diag}((\Sigma(\mu))_{11},1,1)$, $\forall\mu$. Hence, we have from Eq. \eqref{eqn:postDARE} that $\text{trace}(\Sigma^*(\mu_1))=2+\frac{h^2}{\frac{\sigma_1}{2}+h^2}(\sigma_1-1)$,  $\text{trace}(\Sigma^*(\mu_2))=2+\frac{h^2}{\sigma_2+h^2}(\sigma_2-1)$ and $\text{trace}(\Sigma^*(\mu_3))=2+\frac{1}{1-\lambda_1^2}-1$, where $\sigma_1=\sigma(2h^2)$ and $\sigma_2=\sigma(h^2)$ are defined above. Since $\sigma(h^2)$ is a strictly increasing function of $h^2\in\mathbb{R}_{>0}$ with $\sigma(h^2)\ge1$ and upper bounded by $\frac{1}{1-\lambda_1^2}$, and it is easy to obtain that $\frac{\sigma_1}{2}<\sigma_2$, it then follows that Algorithm \ref{algorithm:greedy KFSS} for the posteriori KFSS problem selects sensor $2$ in its first iteration. 

Similarly, we have from Eq. \eqref{eqn:postDARE} that $\text{trace}(\Sigma^*(\mu_{12}))=1+\frac{h^2}{\sigma_2+h^2}(\sigma_2-1)$, $\text{trace}(\Sigma^*(\mu_{23}))=1+\frac{h^2}{2\sigma_{23}+h^2}(\sigma_{23}-1)$ and $\text{trace}(\Sigma^*(\mu_{13}))=1$, where $\sigma_{23}=\sigma(\frac{h^2}{2})$ is defined above. Since $\sigma(h^2)$ is strictly increasing in $h^2\in\mathbb{R}_{>0}$ with $\sigma(h^2)\ge1$ and upper bounded by $\frac{1}{1-\lambda_1^2}$, and it is easy to check that $\sigma_2<2\sigma_{23}$, it follows that the greedy algorithm selects sensor $3$ in its second iteration, and $\mu=\mu_{13}$ is the optimal sensor selection (for the posteriori KFSS instance). Combining the results above and letting $h\to\infty$, we obtain Eq. \eqref{eqn:limit ratio post KFSS}.\hfill\IEEEQED

\subsection*{Proof of Theorem $\ref{thm:ratio of special system KFSA}:$}
We first analyze Algorithm $\ref{algorithm: greedy KFSA}$ for the priori KFSA problem. Since the only nonzero eigenvalue of $A$ is $\lambda_1$, we know from Lemma $\ref{Lemma:minimum trace of sigma}$(a) that $(\Sigma(\nu^c))_{22}=1$ and $(\Sigma(\nu^c))_{33}=1$, $\forall\nu$, which implies that $(\tilde{\Sigma}_{gre})_{22}=1$ and $(\tilde{\Sigma}_{gre})_{33}=1$. Hence, we only need to determine $(\tilde{\Sigma}_{gre})_{11}$.

First, denote $\nu_{1}=[1\ 0\ 0\ 0]^T$, $\nu_{2}=[0\ 1\ 0\ 0]^T$ and $\nu_3=[0\ 0\ 1\ 0]^T$. Then, it is easy to see that $C(\nu^c)$ is of full column rank for all $\nu\in\{\nu_1,\nu_2,\nu_3\}$. This implies that $\mathbf{e}_1\in\text{rowspace}(C(\nu^c))$ for all $\nu\in\{\nu_1,\nu_2,\nu_3\}$. Thus, we know from Lemma $\ref{Lemma:minimum trace of sigma}$(c) that $(\Sigma(\nu^c))_{11}=1$, $\forall\nu\in\{\nu_1,\nu_2,\nu_3\}$. Moreover, denoting $\nu_{4}=[0\ 0\ 0\ 1]^T$, we have $C(\nu_4^c)$ (after some elementary row operations and merging identical rows) is of the form $C(\nu_4^c)=
  \left[\begin{smallmatrix}
    1 & 0 & h \\
    0 & 1 & 0
  \end{smallmatrix}\right]$. Using the results from the proof of Theorem $\ref{thm:ratio of special system KFSS}$, we obtain that $\sigma'_4\triangleq(\Sigma(\nu_4^c))_{11}$ satisfies
\begin{equation*}
\sigma'_4=\frac{1+h^2\lambda_1^2-h^2+\sqrt[]{(h^2-h^2\lambda_1^2-1)^2+4h^2}}{2}.
\end{equation*}
If we view $\sigma'_4$ as a function of $h^2$, denoted as $\sigma'(h^2)$, we know from Lemma $\ref{lemma:estimation error of state 1}$(a) that $\sigma'(h^2)$ is a strictly increasing function of $h^2\in\mathbb{R}_{\ge0}$ with $\sigma'(0)=1$, which implies $\sigma'_4>1$. Thus, Algorithm $\ref{algorithm: greedy KFSA}$ for priori KFSA targets sensor $4$ in its first iteration.

Second, denote $\nu_{14}=[1\ 0\ 0\ 1]^T$, $\nu_{24}=[0\ 1\ 0\ 1]^T$ and $\nu_{34}=[0\ 0\ 1\ 1]^T$. We obtain that $C(\nu^c)$ (after some elementary row operations) is of the form $C(\nu^c)=
  \left[\begin{smallmatrix}
    1 & 0 & h\\
     0 & 1 & 0
  \end{smallmatrix}\right]$, for all $\nu\in\{\nu_{14},\nu_{24},\nu_{34}\}$. It follows that $(\Sigma(\nu^c))_{11}=\sigma'_4$ for all $\nu\in\{\nu_{14},\nu_{24},\nu_{34}\}$, which implies that $\text{trace}(\tilde{\Sigma}_{gre})=\sigma'_4+2$.

Furthermore, the optimal sensor attack (for the priori KFSA instance)  is $\nu=\nu_{12}$, where $\nu_{12}=[1\ 1\ 0\ 0]^T$, since in this case we know from Lemma $\ref{Lemma:minimum trace of sigma}$(a) and (b) that $\Sigma(\nu^{c})_{11}=\frac{1}{1-\lambda_1^2}$, which is also the maximum value of  $\Sigma(\nu^{c})_{11}$ that it can achieve, i.e., $\tilde{\Sigma}_{opt}=\frac{1}{1-\lambda_1^2}+2$. Combining the results above and taking the limit as $h\to 0$, we obtain Eq. $\eqref{eqn:limit ratio KFSA}$.

We next analyze Algorithm \ref{algorithm: greedy KFSA} for the posteriori KFSA problem. Since we know from previous arguments that $C(\nu^c)$ is of full column rank for all $\nu\in\{\nu_1,\nu_2,\nu_3\}$, it follows from Lemma $\ref{Lemma:minimum trace of sigma}$(c) that $\text{trace}(\Sigma^*(\nu_1^c))=\text{trace}(\Sigma^*(\nu_2^c))=\text{trace}(\Sigma^*(\nu_3^c))=0$. Moreover, it is easy to obtain from Eq. \eqref{eqn:couple priori and posteriori 1} that $\Sigma(\nu^c)$ is of the form $\Sigma(\nu^c)=\text{diag}((\Sigma(\nu^c))_{11},1,1)$, $\forall\nu$. We then have from Eq. \eqref{eqn:postDARE} that $\text{trace}(\Sigma^*(\nu_4^c))=1+\frac{h^2}{\sigma'_4+h^2}(\sigma'_4-1)$, where $\sigma'_4=\sigma'(h^2)$ is defined above. Since $\sigma'(h^2)$ is strictly increasing in $h^2\in\mathbb{R}_{\ge0}$ with $\sigma'(0)=1$, it implies that Algorithm \ref{algorithm: greedy KFSA} for posteriori KFSA targets sensor $4$ in its first iteration. Similarly, we have from Eq. \eqref{eqn:postDARE} that $\text{trace}(\Sigma^*(\nu_{14}^c))=\text{trace}(\Sigma^*(\nu_{24}^c))=\text{trace}(\Sigma^*(\nu_{34}^c))=1+\frac{h^2}{\sigma'_4+h^2}(\sigma'_4-1)$, which implies $\text{trace}(\tilde{\Sigma}^*_{gre})=1+\frac{h^2}{\sigma'_4+h^2}(\sigma'_4-1)$.

Furthermore, denote $\nu_{23}=[0\ 1\ 1\ 0]^T$ and $\nu_{13}=[1\ 0\ 1\ 0]^T$. It is easy to show, via similar arguments to those above, that $\text{trace}(\Sigma^*(\nu^c))=1+\frac{h^2}{\sigma'_4+h^2}(\sigma'_4-1)$ for all $\nu\in\{\nu_{34},\nu_{24},\nu_{23},\nu_{14}\}$,  $\text{trace}(\Sigma^*(\nu_{13}^c))=1$, and $\text{trace}(\Sigma^*(\nu_{12}^c))=1+\frac{1}{1-\lambda_1^2}-1$. Since $\sigma'_4=\sigma'(h^2)$ is strictly increasing in $h^2\in\mathbb{R}_{\ge0}$ with $\sigma'(0)=1$, and upper bounded by $\frac{1}{1-\lambda_1^2}$, it follows that the optimal sensor attack (for the posteriori KFSA instance) is $\nu=\nu_{12}$. Combining the results above and taking the limit as $h\to0$, we obtain Eq. \eqref{eqn:limit ratio post KFSA}.\hfill\IEEEQED

\bibliographystyle{IEEEtran}
\bibliography{IEEEabrv,main}
\end{document}